\documentstyle[12pt]{article}

\renewcommand{\H}{\mbox{$\bf H$}}     
\newcommand{\C}{\mbox{$\bf C$}}     
\newcommand{\Z}{\mbox{$\bf Z$}}     
\newcommand{\Q}{\mbox{$\bf Q$}}     
\newcommand{\R}{\mbox{$\bf R$}}     

\newcommand{\cC}{{\cal C}}
\newcommand{\cD}{{\cal D}}
\newcommand{\cE}{{\cal E}}
\newcommand{\cF}{{\cal F}}

\newcommand{\cH}{{\cal H}}
\newcommand{\cI}{{\cal I}}
\newcommand{\cJ}{{\cal J}}
\newcommand{\cK}{{\cal K}}
\newcommand{\cL}{{\cal L}}

\newcommand{\cO}{{\cal O}}

\newcommand{\B}[1]{\begin{#1}}
\newcommand{\E}[1]{\end{#1}}

\newcommand{\eX}{\mbox{$X_{\acute{e}t}$}} 
\newcommand{\Het}[1]{\mbox{$H_{\acute{e}t}^{#1}$}}
\newcommand{\mun}{\mu_n^{\otimes2}}

\newcommand{\im}{{\rm im}\,}        
\newcommand{\coker}{{\rm coker}\,}  
\newcommand{\Pic}{{\rm Pic}\,}      

\newcommand{\Spec}{{\rm Spec}\,}
\newcommand{\df}{\mbox{\,$\stackrel{\pp{\rm def}}{=}$}\,}

\newcommand{\by}[1]{\stackrel{#1}{\rightarrow}}
\newcommand{\longby}[1]{\stackrel{#1}{\longrightarrow}}
\newcommand{\vlongby}[1]{\stackrel{#1}{\mbox{\large{$
\longrightarrow$}}}}

\newcommand{\implies}{\mbox{$\Rightarrow$}}
\newcommand{\tensor}{\otimes}
\newcommand{\into}{\hookrightarrow}

\newcommand{\ie}{{\it i.e.,\/}\ }
\newcommand{\eg}{{\it e.g.\/}\ }
\newcommand{\cf}{{\it cf.\/}\ }

\newcommand{\p}[1]{\mbox{$\scriptstyle {#1}$}}
\newcommand{\pp}[1]{\mbox{$\scriptscriptstyle {#1}$}}

\newcommand{\veq}{\mbox{\large $\parallel$}}

\newcommand{\bul}{\mbox{\Large $\cdot $}}
\newcommand{\Xbul}{X_{\bul}}

\title{Roitman's theorem for singular complex projective surfaces}

\author{by {\sc L.Barbieri-Viale}, {\sc C.Pedrini}${}^*$ and {\sc
C.Weibel}${}^{**}$}

\date{}

\begin{document}

\maketitle

\vspace{2cm}

\B{abstract}
Let $X$ be a complex projective surface with arbitrary
singularities.  We construct a generalized Abel--Jacobi map
$A_0(X)\to J^2(X)$ and show that it is an isomorphism on
{\it torsion\,} subgroups.  Here $A_0(X)$ is the appropriate
Chow group of smooth $0$-cycles of degree 0 on $X$, and
$J^2(X)$ 
is the intermediate Jacobian associated with
the mixed Hodge structure on $H^3(X)$.
Our result generalizes a theorem of Roitman for
smooth surfaces: if $X$ is smooth then the torsion in the
usual Chow group $A_0(X)$ is isomorphic to the torsion
in the usual Albanese variety $J^2(X)\cong Alb(X)$
by the classical Abel-Jacobi map.
\E{abstract}
\vfill
\footnoterule
${}^*$ {\scriptsize Members of GNSAGA of CNR, partially supported by ECC
 Science Plan}

${}^{**}$ {\scriptsize Partially supported by NSF grants}

\newpage

\section*{Introduction}

If $X$ is a smooth projective surface over the
complex numbers $\C$, the classical Abel--Jacobi map
goes from the Chow group $A_0(X)$ of cycles of degree 0
to the (group underlying the) Albanese Variety $Alb(X)$.
Roitman's Theorem \cite{Roit} states that this map induces
an isomorphism on torsion subgroups.
(See \cite{CT} for a nice compendium).

The goal of this paper is to remove the word ``smooth''
 from Roitman's theorem.  For this we shall modify
the definition of $A_0(X)$, replace $Alb(X)$ with Griffiths'
intermediate Jacobian $J^2(X)$, and construct a
generalization of the Abel--Jacobi map.
\medskip

\noindent{\bf Main Theorem. }\
{\it Let $X$ be a reduced projective
surface over $\C$. Then there is a natural map from $A_0(X)$
to $J^2(X)$ inducing an isomorphism on torsion:
$$A_0(X)_{tors}\cong J^2(X)_{tors}.$$
In particular, the torsion subgroup is a finite direct sum
of copies of $\Q/\Z$.}

\smallskip
If $X$ is a normal surface, this theorem is a reformulation of
a theorem of Collino  and Levine \cite{C2} \cite{L-Alb}, because
(as we will show in Corollary~\ref{abnorm}),  $J^2(X)$
is isomorphic to the Albanese of any desingularization of $X$.

Gillet studied the Abel--Jacobi map in \cite{GDuke} when
$X$ is a singular surface with ``ordinary multiple curves''
(\eg a seminormal surface with smooth normalization $\tilde X$).
He proved in \cite[Theorem B]{GDuke} that if $\tilde X$ satisfied
some extra hypotheses ($p_g=0$, etc.) then the Abel--Jacobi map
is surjective with finite kernel.  Thus we deduce:

\medskip
\noindent{\bf Corollary. }\  {\it
Let $X$ be a surface with ordinary multiple curves
such that $H^2(X,\cO_X)=0$.  Assume that Bloch's conjecture
holds for the normalization $\tilde X$ of $X$.
Then the Abel--Jacobi map is an isomorphism.
$$A_0(X) \cong J^2(X)$$}

We now describe the ingredients in our main theorem.
If $X$ is a proper surface over $\C$, the intermediate
Jacobian $J^2(X)$ is defined to be
$$J^2(X) \df\frac{H^{3}(X,\C)}{F^2H^{3}+ H^{3}(X,\Z(2))}.$$
Here $F^2H^3$ refers to the Hodge filtration of \cite{D}
and the coefficients $\Z(2)$ refer to the
embedding of $\Z$ in $\C$ sending 1 to $(2\pi i)^2$.

When $X$ is a smooth surface, it is well known that
$J^2(X)$ is isomorphic to the Albanese Variety $Alb(X)$.
Now suppose that $X$ is a singular surface.
We will show in Corollary~\ref{1-motive} that $J^2(X)$
is a complex torus, and that
if $X'$ is a resolution of singularities for $X$
then $J^2(X)$ is an extension of $Alb(X')$ by a torus.
That is, the map $J^2(X)\to Alb(X')$ forms a 1--motive
in the sense of Deligne \cite{D}; we call it the
{\it Albanese 1--motive} of $X$.
Given this, the torsion subgroup of $J^2(X)$
is a finite direct sum of copies of $\Q/\Z$.

The modified version of $A_0(X)$ is defined as a subgroup
of the Levine--Weibel Chow group $CH_0(X)$ of zero-cycles
on $X$ \cite{LW}.  By definition, $CH_0(X)$ is the
abelian group generated by the smooth closed points on $X$,
modulo the subgroup generated by all terms $D=\sum n_iP_i$
(with $P_i$ smooth on $X$) such that
$D = (f)$ for some rational function $f$ on some
curve $C$, the curve being locally defined by a single equation
on the surface $X$.

If $X$ is a surface with $c$ proper components, there is a
natural surjection $CH_0(X)\to\Z^c$, called the degree map.
By definition, $A_0(X)$ is the kernel of the degree map.

In order to prove our Main Theorem, we need to reinterpret
$CH_0(X)$ in terms of algebraic $K$-theory.
Let $SK_0(X)$ denote the subgroup of $K_0(X)$
generated by the classes of smooth points on the surface $X$.
Then $CH_0(X)$ is isomorphic to $SK_0(X)$, by the map sending
a smooth point to its class in $K_0(X)$.  This is the
Riemann--Roch Theorem if $X$ is smooth.  It was proven for
affine surfaces in \cite[Theorem 2.3]{LW}.
For arbitrary quasi-projective surfaces it is due to
Levine \cite{L1}, who proved that both groups are isomorphic
to $H^2(X,\cK_2)$ (\cf \cite{PW1}, \cite{C1}).
The isomorphism
$$CH_0(X) \cong H^2(X,\cK_2) \cong SK_0(X)$$
is often called ``Bloch's formula'' for surfaces.

We have laid this paper out as follows.
In \S1 we present some basic facts about Deligne cohomology
of a proper but singular scheme.  The corresponding Deligne
Chern classes which will be used in later sections is
introduced in \S2.
In \S3 we construct and compare the Mayer-Vietoris sequences
for $K$-theory and Deligne cohomology that we shall need.
In \S4 we compute $J^2(X)$ for any proper surface $X$.
Our computation shows that $J^2(X)$ is part of a 1-motive
$Alb(X)$ which we call the {\it Albanese 1-motive of $X$}.
In \S5 we describe the structure of $SK_1$ of any curve over
any algebraically closed field.
In \S6 and \S7 we establish some technical results about
$\cK_2$-cohomology, ending with the exact sequence of
Theorem~\ref{NH3}
for a normal surface $X$ over any field containing $\frac1n$.
$$0\to H^1(X,\cK_2)/n\to N\Het3(X)\to{}_nCH_0(X)\to 0.$$
Finally we prove the Main Theorem in \S8.

\medskip
\section*{Notation}

All schemes we consider will be separated and
of finite type over a field $k$.
We call such a scheme a {\it curve} if it is 1-dimensional,
and a {\it surface} if it is 2-dimensional.
When $X$ is an algebraic scheme over $\C$, we will write
$H^*(X,\Z)$ and $H^*(X,\C)$ for the singular cohomology of the
associated analytic space $X_{an}$ as well as for
the mixed Hodge structure on it, given by Deligne \cite{D}.
The weight filtration on $H^*(X,\Z)$ will be written as
$W_iH^*$, and the Hodge filtration on $H^*(X,\C)$
will be written as $F^iH^*$.

The notation $\Z(r)$ denotes the subgroup $(2\pi i)^r\Z$ of $\C$.
Unless we wish to call attention to the relation with $H^*(X,\C)$,
we will write $H^*(X,\Z)$ instead of $H^*(X,\Z(r))$.
The notation $\Z(r)_\cD$ denotes the Deligne complex
on a smooth scheme $X$ over $\C$ (see \S1).
We will use the Deligne complex to define the Deligne cohomology of
proper schemes; in the affine case the definition of
Deligne-Beilinson cohomology is different (one needs to consider
logarithmic poles), and we remain silent about this.
Similarly, the Zariski sheaves
$\cH_\cD^*(r)$ (defined as the higher direct images of $\Z(r)_\cD$)
are used only for proper schemes, as a technical device.
(See \S1, (2.4), \ref{square} and \ref{crux}.)

The Zariski sheaf $\cK_q$ on $X$ is obtained by sheafifying
the Quillen higher $K$-theory functor $U\mapsto K_q(U)$. The
$\cK$-cohomology groups $H^p(X,\cK_q)$ are just the
Zariski cohomology of these sheaves.   As indicated in the
introduction, when $X$ is a surface the most important
$\cK$-cohomology group is $H^2(X,\cK_2)\cong CH_0(X)$.

Similarly, we shall write $\cK_q(\Z/n)$ and
$\cH^q(\mu_n^{\otimes i})$ for the Zariski sheaves associated
to the presheaves sending $U$ to $K_q(U;\Z/n)$ and
$\Het{q}(U,\mu_n^{\otimes i})$, respectively.  In general,
we will always use calligraphic letters for Zariski sheaves.

Finally, we will use some standard notation.
Let $H$ be an abelian group or sheaf of abelian groups.
Then $H_{tors}$ will denote its torsion subgroup.
For each integer $n$ we will write $H/n$ for $H/nH$,
and ${}_nH$ for the subgroup $\{ x\in H\colon nx=0\}$ of $H$.

\goodbreak
\section{Deligne cohomology groups}

For $X$ smooth (possibly affine) over $\C$ we let
$\Z(r)_{\cD}$ denote the ``Deligne complex''
$$0\to \Z(r)\to \cO_{X_{an}}\by{d} \cdots \by{d}
\Omega^{r-1}_{X_{an}} \to0$$
of sheaves on the complex analytic manifold $X_{an}$,
where $\Z(r)\df (2\pi i)^r\Z$ is in degree 0.
The {\it analytic} Deligne cohomology groups of the smooth
scheme $X$ are defined to be
$$H_\cD^q(X,\Z(r))=H^q(X,\Z(r)_{\cD})
	\df\H_{an}^q(X,\Z(r)_{\cD}).$$
We then have exact sequences of complexes of sheaves
on $X_{an}$:
\B{equation}\label{augment}
0\to \Omega^{<r}_{X_{an}}[-1]\to
\Z(r)_{\cD}\by{\varepsilon_X} \Z\to 0.
\E{equation}
We can also define the Deligne cohomology groups of a smooth
simplicial scheme $\Xbul$ by considering
$\Z(r)_{\cD}$ as a complex of analytic sheaves on $\Xbul$.
This yields an exact sequence of complexes parallel to
(\ref{augment}) by \cite[5.1.9.(II)]{D}.

Now let $X$ be a singular scheme. 
A {\it smooth proper hypercovering} $\Xbul\to X$ of
$X$ (\cf \cite[6.2.5--6.2.8]{D}) is a simplicial scheme
$\Xbul$ with smooth components $X_i$, each proper over $X$,
together with a morphism to $X$ satisfying ``universal
cohomological descent.''
We define the Deligne cohomology of $X$ to be:
$$H_{\cD}^q(X,\Z(r))\df\H_{an}^q(\Xbul,\Z(r)_{\cD}).$$
This definition is independent of the choice of smooth proper
hypercovering by \cite[Expos\'e V{\it bis},
5.1.7 and 5.2.4]{SGA4}.
There is a canonical descent isomorphism
$H^*(X,\Z)\cong H^*(\Xbul,\Z)$, so the map $\varepsilon$ in
(\ref{augment}) induces a natural map
$\varepsilon_X\colon H_\cD^*(X,\Z(r))\to H_{an}^*(X,\Z)$.
It is well-known (see \cite[1.6.4]{Bei}{})
that $\varepsilon_X$ preserves products.

For $X$ proper with arbitrary singularities
we have a standard long exact sequence
\B{equation}\label{modf}
\cdots\by{\varepsilon}\kern-2pt H^q(X,\Z)\to\kern-2pt
H^q(X,\C)/F^r\kern-2pt\to\kern-2pt H_{\cD}^{q+1}(X,\Z(r))
\by{\varepsilon}\kern-2pt H^{q+1}(X,\Z)\to\cdots
\E{equation}
induced by (\ref{augment}) and $\Z\cong\Z(r)\subset\C$,
as well as
\B{equation}\label{star}
\cdots\to F^rH^q(X,\C)\to H^q(X,\C/\Z(r)) \to
H_{\cD}^{q+1}(X,\Z(r)) \to F^rH^{q+1}(X,\C)\to\cdots
\E{equation}
If $X$ is a proper surface then from (\ref{modf})
we have an exact sequence
\B{equation}\label{extJZ}
0\to J^2(X) \to H_{\cD}^{4}(X,\Z(2))
\by{\varepsilon} H^{4}(X,\Z(2))\to0.
\E{equation}

Any map $i\colon Y\to X$ lifts to a morphism
$i\colon Y_{\bul}\to \Xbul$ between hypercoverings;
see \cite[Expos\'e V{\it bis}, 5.1.7 and 5.2.4]{SGA4}
or \cite[6.2.8]{D}.
The {\it relative} Deligne cohomology of this map is
defined in the notation of \cite[6.3.3]{D} to be:
$$H_{\cD}^q(X{\rm mod\,}Y,\Z(r))\df
\H_{an}^q(\Xbul{\rm\,mod\,}Y_{\bul},\ \Z(r)_{\cD}
{\rm\, mod\,}\Z(r)_{\cD})
$$
By \cite[6.3.2.2]{D} we have a functorial long exact sequence
\B{equation} \cdots\to  H_{\cD}^q(X{\rm\,mod\,}Y,\Z(r))\to
H_{\cD}^q(X,\Z(r))\to H_{\cD}^q(Y,\Z(r)) \to\cdots
\E{equation}
and of course there are relative versions
of (\ref{modf}) and (\ref{star})
which depend functorially on the pair $(X,Y)$, such as
\B{equation}\label{modf-rel}
\cdots\by\varepsilon\kern-3pt H^q(X{\rm mod\,}Y,\Z)\to\kern-3pt
H^q(X{\rm mod\,}Y )/F^r\to\kern-3pt
H_{\cD}^{q+1}(X{\rm mod\,}Y,\Z(r))\by\varepsilon\cdots.
\E{equation}
\subsection*{Low degree Deligne cohomology} 

We will need the following calculation of $H_\cD^q(X,\Z(2))$
for $q\le2$.  Given a proper scheme $X$ over $\C$, we fix a smooth
proper hypercovering $\Xbul\by\pi X$.  By abuse of notation,
we shall write $\cH_\cD^q(2)$ for the complexes $R^q\omega_*\Z(2)_\cD$
of Zariski sheaves on either $\Xbul$ or $X$,
$\omega$ denoting either
$\omega_{\bul}\colon X_{\bul an}\to X_{\bul zar}$ or
$\omega=\pi\omega_{\bul}\colon X_{\bul an}\to X_{zar}$.

\B{prop}\label{sheafHD}
For $X$ proper and connected over $\C$ we have:

\kern-8pt
\B{description}
\item[{\it (i)}]  $H_\cD^0(X,\Z(2)) = 0$;
\item[{\it (ii)}]  $H_\cD^1(X,\Z(2))\cong\C/\Z(2) \cong\C^*$;
\item[{\it (iii)}] $H_\cD^2(X,\Z(2))_{tors}\cong H^1(X,\Q/\Z)$.
\E{description}\kern-8pt
Moreover, if $X$ is irreducible then we have
\B{description}
\item[{\it (iv)}]  $H_\cD^2(X,\Z(2)) = H^0_{zar}(X,\cH_\cD^2(2))
=H^0_{zar}(\Xbul,\cH_\cD^2(2))$
\E{description}
and there are edge homomorphisms:
\B{description}
\item[{\it (v)}] $H^1_{zar}(X,\cH_\cD^2(2)) \into
H^1_{Zar}(\Xbul,\cH_\cD^2(2))\into H_\cD^3(X,\Z(2))$
\quad (these are injections);
\item[{\it (vi)}]  $H^2_{zar}(X,\cH_\cD^2(2)) \longby{}
H^2_{Zar}(\Xbul,\cH_\cD^2(2))\to H_\cD^4(X,\Z(2))$.
\E{description}
\E{prop}

\B{proof}It is well-known that $H_{an}^1(X,\Z)$ is torsion-free.
Hence $(i)$ and $(ii)$ follow immediately from (\ref{modf}).
(Cf. the proof of Lemma~2.17 in \cite{GDuke}.)
Part $(iii)$ follows from this and (\ref{star}),
since $H^1(X,\C/\Z(2))_{tors}\cong H^1(X,\Q/\Z)$.
In order to prove parts $(iv)$, $(v)$ and $(vi)$
 we use the Leray spectral sequences
for $\omega$ and $\omega_{\bul}$:
\B{equation} \label{simps}
\B{array}{ccccc}
E^{p,q}_2 &=& H^p_{zar}(X,\,\cH_\cD^q(i))
&\implies& H_\cD^{p+q}(X,\Z(i))_{\strut} \\
E^{p,q}_2 &=& H^p_{zar}(\Xbul,\cH_\cD^q(i))
&\implies& H_\cD^{p+q}(X,\Z(i))
\E{array}\E{equation}
with $i=2$ (\cf \cite[(2.13)]{GDuke}).
For this, we need to compute $\cH_\cD^0(2)$ and $\cH_\cD^1(2)$.

 When $U$ is smooth we may identify the analytic sheaf
$\cO_U/\Z(2)$ with $\cO_U^*$ and
obtain a quasi-isomorphism between $\Z(2)_\cD$ and the complex
$0\to\cO_U^*\longby{dlog}\Omega_U^1$.
It follows that there is a distinguished triangle
of complexes of analytic sheaves on $\Xbul$
$$\C^*[-1] \to \Z(2)_\cD \to
\Omega_{\Xbul}^1/dlog(\cO_{\Xbul}^*)[-2]
\to \C^*.$$
Applying $\omega_*$ and $R^1\omega_*$ immediately yields
$\cH_\cD^0(2)=0$ and $\cH_\cD^1(2)=\omega_*\C^*=\C^*$ on both
$X_{\bul zar}$ and $X_{zar}$.  Therefore
in either spectral sequence (\ref{simps})
the row $q=0$ vanishes and in row $q=1$ we have
$H_{zar}^p(X,\C^*)=H_{zar}^p(\Xbul,\C^*)$.
The exact sequences of low degree terms in (\ref{simps}) become:
$$\B{array}{c}
0\to H_{zar}^1(X,\C^*)\to H_\cD^2(X,\Z(2)) \to
 H^0(X,\cH_\cD^2) \longby{d_2} H_{zar}^2(X,\C^*)
\to H_\cD^3(X,\Z(2))\\
0\to H_{zar}^1(X,\C^*)\to H_\cD^2(X,\Z(2)) \to
H^0(\Xbul,\cH_\cD^2) \longby{d_2} H_{zar}^2(X,\C^*)
\to H_\cD^3(X,\Z(2)).
\E{array}$$
The map between these sequences identifies them, and
$H^0(X,\cH_\cD^2)\cong H^0(\Xbul,\cH_\cD^2)$ by the 5-lemma.
If $X$ is irreducible then $H_{zar}^p(X,\C^*)=0$ for $p\ne0$.
Hence $H_\cD^2(X,\Z(2))$ is isomorphic to
$H^0(X,\cH_\cD^2)$.
Parts $(v)$ and $(vi)$ follows similarly.
 \E{proof}

For each $n$ there is a distinguished triangle of
complexes of analytic sheaves on $\Xbul$:
\B{equation}\label{triangle}
\Z/n[-1]\by\delta\Z(2)_\cD\by{n}\Z(2)_\cD
\by{\bar\varepsilon}\Z/n \longby{\delta[1]}\Z(2)_\cD[1].
\E{equation}
The comparison theorem between the analytic and \'etale sites,
together with universal cohomological descent, yields
$H^q(\Xbul,\Z/n)\cong H^q(X,\Z/n)\cong\Het{q}(X,\Z/n)$.
Fixing an $n^{th}$ root of unity allows us to identify
$\mu_n$, $\mun$ and $\Z/n$ on $\eX$.
If $X$ is proper, the cohomology of the triangle
(\ref{triangle}) yields ``Kummer sequences''
\B{equation}\label{Kum-etHD}
0\to H_\cD^{q}(X,\Z(2))/n \longby{\bar\varepsilon}
\Het{q}(X,\mun)\longby\delta H_\cD^{q+1}(X,\Z(2))_{n-tors} \to0.
\E{equation}
By \ref{sheafHD} this identifies $\mu_n\cong\Het0(X,\Z/n)$
with the $n$-torsion in $\C^*\cong H_\cD^1(X,\Z(2))$,
and identifies $\Het1(X,\mun)\cong\Het1(X,\Z/n)$
with the $n$-torsion in $H_\cD^2(X,\Z(2))$.

\smallskip
Now consider the morphism $\omega\colon X_{\bul an}\to X_{zar}$.
Applying the higher direct image $R^q\omega_*$ to
(\ref{triangle}) yields an exact sequence of Zariski sheaves:
\B{equation}\label{Rq}
\longby{\bar\varepsilon} \cH^{q-1}(\mun)
\longby\delta \cH_\cD^q(2) \longby{n} \cH_\cD^q(2)
\longby{\bar\varepsilon} \cH^q(\mun)\longby\delta.
\E{equation}
In particular, $\delta$ identifies $\cH^1(\mun)$ with the
$n$-torsion subsheaf of $\cH_\cD^2(2)$.

The map $\Het{*-1}(X,\mun)\by\delta H_\cD^*(X,\Z(2))$
is also the abutment of a morphism of Leray spectral sequences.
At the $E_2$-level it is
$H_{zar}^p(X,\cH^{q-1}(\mun)) \by\delta H_{zar}^p(X,\cH_\cD^q(2))$.
If $X$ is proper and irreducible then the bottom row of both
spectral sequences degenerates
(\eg $H_{zar}^p(X,\Z/n)=0$ for $p\ne0$)
and we obtain the following result.

\B{cor}\label{morph}If $X$ is proper and irreducible,
there is a commutative diagram whose rows are the
exact sequences of low degree terms of Leray spectral sequences:
$$\B{array}{ccccccccc}
0&\to& \kern-5pt H^1(X,\cH^1(\mun))&\to&\kern-3pt\Het2(X,\mun)
&\to& \kern-3pt H^0(X,\cH^2(\mun))
&\kern-4pt\longby{d_2}&\kern-8pt H^2(X,\cH^1(\mun))\\
&&\delta\downarrow &&\delta\downarrow&&\delta\downarrow&&
\delta\downarrow\\
0&\to& \kern-5pt H^1(X,\cH_\cD^2(2))&\to&\kern-3pt
H_\cD^3(X,\Z(2)) &\to& \kern-3pt H^0(X,\cH_\cD^3(2))
&\kern-4pt\longby{d_2}& \kern-8pt H^2(X,\cH_\cD^2(2)).
\E{array}$$
\E{cor}

\section{Chern classes in Deligne cohomology}

For each scheme $X$ of finite type over $\C$ the
exponential map $\cO_{X_{an}}\to\cO_{X_{an}}^*$ induces a
quasi-isomorphism between $\Z(1)_{\cD}=(\Z\to\cO_{X_{an}})$ and
$\cO_{X_{an}}^*[-1]$.  This quasi-isomorphism also holds over a
simplicial scheme $\Xbul$ by naturality, so
$\H^q(\Xbul,\Z(1)_\cD)\cong H^{q-1}(\Xbul,\cO^*_{\Xbul})$.
This gives a natural map from $H^1_{an}(X,\cO_X^*)$ to
$\H^2(\Xbul,\Z(1)_\cD)\cong H_{an}^1(\Xbul,\cO^*_{\Xbul})$
for every smooth proper hypercovering $\Xbul\to X$.
Composing with the determinant map $K_0(X)\to \mbox{Pic}(X)$
and the natural map $\mbox{Pic}(X)\to H_{an}^1(X,\cO_X^*)$
yields a map $c_1\colon K_0(X) \to \H^2(\Xbul,\Z(1)_\cD)$.

Now the splitting principle holds for Deligne cohomology by
\cite[5.2]{GDuke}.  (Warning: if $X$ is not proper this differs
slightly from the splitting principle proven in \cite[1.7.2]{Bei}!)
Thus the map $c_1$ extends to Chern classes
$c_i\colon K_0(X) \to \H^{2i}(\Xbul,\Z(i)_\cD)$
for vector bundles.
When $X$ is proper, these are the Deligne-Beilinson Chern classes
$$c_i\colon K_0(X) \to H_\cD^{2i}(X,\Z(i)).$$
Recall from (\ref{augment}) that there is a map
$\varepsilon_X\colon \H^{2i}(\Xbul,\Z(i)_\cD)\to
H_{an}^{2i}(X,\Z)$, and that it is product-preserving.

\B{lemma}\label{c2-an} (\cf Beilinson \cite[1.7]{Bei})
The composition of $c_i$ with the map $\varepsilon_X$
is the classical Chern class
of the associated topological vector bundle \cite{MCC}:
$$c_i^{an}\colon
	K_0(X)\to H_{an}^{2i}(X,\Z)=H_{top}^{2i}(X,\Z).$$
\E{lemma}
\B{proof} Since $\varepsilon_X$ preserves cup products, the
splitting principle shows that it suffices to establish the
result for $c_1$.
If $X$ is smooth then $c_1$ is the analytic determinant map,
and $\varepsilon_X$ is just the usual map
$\partial_X\colon H_{an}^1(X,\cO_X^*)\to H_{an}^2(X,\Z)$
used to define $c_1^{an}$ on analytic vector bundles,
so it is clear that $c_1^{an}=\varepsilon_X\circ c_1$.
To deduce the result for general $X$, choose a smooth proper
hypercover $u\colon \Xbul\to X$.  Composing $\partial_X$
(which is $c_1^{an}$) with the descent isomorphism
$H_{an}^{2i}(X,\Z)\cong H_{an}^{2i}(\Xbul,\Z)$ is the
descent map
$H^1(X,\cO_X^*)\to H^1(\Xbul,\cO_{\Xbul}^*)$
(which is $c_1$) composed with $\partial_{\Xbul}$,
\ie with $\varepsilon_{\Xbul}$.
\E{proof}
Reduction of $\varepsilon_X$ mod $n$ yields a map
$\bar\varepsilon_X\colon H_\cD^{2i}(X,\Z(i))\to
H_{an}^{2i}(X,\Z/n)$.  Since reduction mod $n$ is
product-preserving and sends $c_1^{an}$ to the
\'etale Chern class $c_1^{et}$, we deduce the
\B{cor}\label{c-etale}
The composition of $c_i$ with $\bar\varepsilon_X$
is the \'etale Chern class
$$c_i^{et}\colon K_0(X) \to H_{an}^{2i}(X,\Z/n)
\cong\Het{2i}(X,\mu_n^{\otimes i}).$$
\E{cor}

In this paper we shall be mostly concerned with the class
$c_2\colon K_0(X)\to H^4_{\cD}(X,\Z(2))$ when $X$ is a
projective surface.  Recall from the introduction (or \cite{LW})
that the Chow group $CH_0(X)$ of zero-cycles on $X$ is
isomorphic to the subgroup $SK_0(X)$ of $K_0(X)$.

If $X$ has $c$ irreducible components then there is a natural
degree map $CH_0(X)\to\Z^c$, and $A_0(X)$ is defined to be
the kernel of this map.
The following cohomological interpretation
of the degree map will be useful.

\B{lemma}\label{degree} (Beilinson \cite[1.9]{Bei})
If $X$ is a projective surface, the
degree map is the same (up to sign) as the classical Chern class
$$ CH_0(X) \hookrightarrow K_0(X) \by{c_2^{an}} H_{an}^4(X,\Z)
	\cong \Z^c.$$
By (\ref{extJZ}), the Deligne Chern class $c_2$ induces a natural map
$\rho\colon A_0(X)\to J^2(X)$ 
fitting into the diagram
$$ \begin{array}{ccccccc}
0 \to & A_0(X) &\to& CH_0(X) &\by{deg}& \Z^c &\to0 \\
  &\rho\downarrow& &c_2\downarrow&    & \veq &     \\
0 \to & J^2(X) &\to& H_\cD^4(X,\Z(2)) &\by{\varepsilon}&
H^4_{an}(X,\Z) &\to0 \\
\end{array}$$
\E{lemma}

\smallskip\noindent{\bf Definition: }
We shall refer to the map $\rho$ as the {\it Abel--Jacobi map},
because if $X$ is a smooth surface then
$J^2(X)$ is the usual Albanese variety and the map
$\rho$ coincides with the classical Abel--Jacobi map
by \cite[1.9.1]{Bei} or \cite[2.24]{GDuke}.

\medskip
\B{proof} Observe that if $X$ has $c$ proper irreducible
components then $H^4(X,\Z)\cong\Z^c$, because
the singular locus of $X$ has real analytic dimension $\le2$.
Given Lemma~\ref{c2-an}, the second assertion follows from
the first.  If $X$ is a smooth projective surface the result is
classical; one way to see it is to use the product formula
for two divisors on $X$:
$$c_2^{an}(D\otimes E) = -c_1^{an}(D)\cup c_1^{an}(E)
			= - (D\cdot E)[X].$$
In general, choose a resolution of singularities $X'\to X$.
Since $X'$ has $c$ disjoint components, the degree map on
$X$ factors through the degree map on $X'$ as
$CH_0(X)\to CH_0(\tilde X)\to\Z^c$.
By naturality of $c_2^{an}$, the isomorphism
$H^4_{an}(X,\Z)\cong H^4_{an}(X',\Z)$ allows us to
deduce the result for $X$ from the result for $X'$.
\E{proof}

\B{num}\label{simplicial}
As observed by Beilinson \cite[2.3]{Bei}
(\cf \cite[\S5]{GDuke}), the formalism of
Deligne cohomology allows us to extend the Chern classes
 from $K_0(X)$ to higher $K$-theory as well.  The higher
Deligne Chern classes are homomorphisms
$$c_i\colon K_q(X) \to H^{2i-q}_\cD(X,\Z(i)).$$
Composition with $\varepsilon_X$ yields the higher analytic
Chern classes $c^{an}_i$, and reduction mod $n$ yields
the higher \'etale Chern classes $c^{et}_i$.
Moreover, the following holds.

\B{description}
\item{(\thethm.1)}
There is a connected simplicial presheaf
$K\simeq\Omega_0BQP$ and a simplicial sheaf $\cD$
on $X_{zar}$ such that $\pi_q K(U) = K_q(U)$ for $q\ge1$,
and $\pi_q\cD(U) = H^{2i-q}(U_{\bul},\Z(i)_\cD)$ for $q\ge0$.
Moreover, there is a map of simplicial presheaves
$C_i^{ss}\colon  K \to \cD$
such that  $\pi_q C_i^{ss}(X)$ is the
Deligne cohomology Chern class $c_i$ on $K_q(X)$.
(\cf \cite[5.4]{GDuke}, which differs somewhat from
\cite{Bei} and \cite{GBAMS}.)

Indeed, $\cD$ is the simplicial sheaf of abelian groups associated
by the Dold-Kan theorem (\cite[8.4.1]{W-homo})
to the good truncation  $\tau^{\le0}\R\omega_*\Z(i)_\cD[2i]$
of the total derived direct image of $\Z(i)_\cD[2i]$ under
$\omega\colon {\Xbul}{}_{,an}\to X_{zar}$.

\item{(\thethm.2)} Let $\cE$ denote the simplicial sheaf associated
by the Dold-Kan theorem to the good truncation
$\tau^{\le0}\R\omega_*\Z/n[2i]$ of the total derived direct
image of $\Z/n[2i]$.  Then $\pi_q\cE(U)=H_{an}^{2i-q}(U,\Z/n)
\cong\Het{2i-q}(U,\mu_n^{\otimes i})$.  If we define $L$ to be
the homotopy fiber of $K\by{n}K$ then we have
$\pi_q L(U) = K_{q+1}(U;\Z/n)$.  This all gives
a homotopy commutative diagram whose rows are
homotopy fibration sequences
\B{equation}\label{fib}\B{array}{ccccccccc}
\Omega K &\to& L &\to& K &\by{n}& K && \\
\kern18pt\downarrow\Omega C^{ss}_2\kern-10pt &&
\kern13pt\downarrow C^{ss}_2 \kern-10pt&&
\kern13pt\downarrow C^{ss}_2 \kern-10pt&&
\kern13pt\downarrow C^{ss}_2 \kern-8pt&& \\
\Omega\cD &\to& \Omega\cE &\by\delta&\cD&\by{n}&\cD&
\by{\bar\varepsilon}&\cE.\\
\E{array}\E{equation}
 From Corollary~\ref{c-etale}
and a standard argument with $\Het{*}(X,G,\mu_n^{\otimes i})$
it is easy to see that not only does
$K\to\cE$ induce the higher \'etale Chern class $c^{et}_i$
on $K_*(X)$ but the map $L\to\Omega\cE$ induces
the usual \'etale Chern classes on
$K$-theory with coefficients mod $n$.
$$c^{et}_i\colon K_q(X;\Z/n)\longby{}
\Het{2i-q}(X,\mu_n^{\otimes i})
$$
Applying $\pi_2$ to (\ref{fib}) with $i=2$ and $U=X$ yields
the commutative diagram
\B{equation}\label{et-HD}\B{array}{ccccccc}
K_3(X)&\to&K_3(X;\Z/n) &\longby{}& K_2(X) &\by{n}& K_2(X) \\
&&\downarrow c^{et}_2 && \downarrow c_2 && \downarrow c_2 \\
0&\to& \Het1(X,\mun) &\longby{\delta}& H_\cD^2(X,\Z(2))
&\by{n}& H_\cD^2(X,\Z(2)).
\E{array}\E{equation}
By (\ref{Kum-etHD}) we see that $c^{et}_2$ vanishes on
$K_3(X)$ and factors through ${}_nK_2(X)$.

Applying $\pi_2$ to (\ref{fib}) with $i=2$ and
sheafifying yields the commutative diagram of sheaves
in which the bottom row is part of (\ref{Rq}):

\B{equation}\label{et-HD-sheaf}\B{array}{ccccccc}
\cK_3&\to&\cK_3(\Z/n) &\longby{}&
\kern-5pt\cK_2 &\by{n}& \kern-5pt\cK_2 \\
&&\downarrow c^{et}_2 && \downarrow c_2 && \downarrow c_2 \\
0&\to& \cH^1(\mun) &\longby{\delta}& \cH_\cD^2(2)
&\by{n}& \cH_\cD^2(2)
\E{array}\E{equation}
By (\ref{Rq}) we see that $c^{et}_2$ vanishes on
$\cK_3$ and factors through the sheaf ${}_n\cK_2$.

\item{(\thethm.3)}
There is a morphism of spectral sequences
between the Brown--Gersten spectral sequence for $K_{-*}(X)$
and the Leray spectral sequence in (\ref{simps}) converging to
 $H^{2i+*}_\cD(X,\Z(i))$.  At the $E_2^{pq}$-level the
morphisms are the cohomology of $c_i$:
$$H^p_{zar}(X,\cK_{-q}) \by{c_i} H^p_{zar}(X,\cH_\cD^{2i+q}(i))
.$$
Here $\cK_q$ is the sheaf on $X_{zar}$
associated to the presheaf $K_q$ and the sheaves
$\cH_\cD^j(i)$ are $R^j\omega_*\Z(i)_\cD$, as in the proof
of Proposition~\ref{sheafHD}.
By \cite{TT}, the first spectral sequence converges to
$K_{-p-q}(X)$ whenever $X$ is quasi-projective.
The second spectral sequence is an obvious
reindexing of (\ref{simps}) and converges to
$H_\cD^{2i+p+q}(X,\Z(i))$.
\E{description}
\E{num}
Here are three applications of the morphism of spectral sequences
in (2.4.3).  First, if $X$ is a projective surface
we have a commutative diagram
$$\begin{array}{ccccc}
CH_0(X)&\cong\kern-8pt& H^2(X,\cK_2) &\hookrightarrow& K_0(X)\\
\strut     & & c_2\downarrow&         & c_2\downarrow   \\
      & & H^2(X,\cH_\cD^2(2)) &\to  & H_\cD^4(X,\Z(2)).\\
\end{array}$$
where the bottom horizontal map is given by
Proposition~\ref{sheafHD}$(vi)$.

Second, suppose that $Y$ is 1-dimensional. \kern-1pt Then we may
identify the group $H^1\kern-1pt(X\kern-1pt,\cK_2\kern-1pt)$
with the subgroup $SK_1(X)$ of $K_1(X)$,
and $c_2\colon SK_1(X)\to H_\cD^3(X,\Z(2))$ is identified with
the composite
$H^1(X,\cK_2)\by{c_2}H^1(X,\cH_\cD^2(2))\to H_\cD^3(X,\Z(2))$.

Third, suppose that $X$ is an irreducible projective surface.
Then $c_2$ vanishes on the image of $H^2(X,\cK_3)$ in $K_1(X)$
because it factors through $H^2(X,\cH_\cD^1(2))=H_{zar}^2(X,\C^*)$,
which is zero because $X$ is irreducible, as we saw in the
proof of Proposition~\ref{sheafHD}.
Since $SK_1(X)$ is an extension of $H^1(X,\cK_2)$ by this image,
we may summarize this as follows.

\B{lemma}\label{c2-SK1}
Let $X$ be an irreducible projective surface over $\C$.  Then
the Chern class
$c_2\colon SK_1(X)\to H_\cD^3(X,\Z(2))$ factors as:
$$SK_1(X)\vlongby{\mbox{\rm onto}}H^1(X,\cK_2)
\by{c_2}H^1(X,\cH_\cD^2(2))\into H_\cD^3(X,\Z(2)).
$$
\E{lemma}

\goodbreak
\section{Mayer--Vietoris sequences} 

Since we are going to deal with resolutions of singularities or
normalizations we will need some Mayer--Vietoris sequences.
In this section we do this for mixed Hodge structures,
Deligne cohomology and $K$-theory.

Associated to a proper birational morphism
$f\colon X'\to X$ of $\C$-algebraic schemes, and every
closed subscheme $i\colon Y\into X$ we have the  commutative square
\B{equation}\label{birsquare}
\B{array}{rcl} Y'& {\stackrel{i'}{\into}}& X'\\
\p{f'}\downarrow & &\downarrow\p{f}\\
Y &{\stackrel{i}{\into}}& X  \E{array}
\E{equation}
where $Y' = f^{-1}(Y)\ (\ =Y\times_X X')$.
We shall always assume that $Y$ is chosen so that the restriction
$f\colon X'-Y'\by{\simeq} X-Y$ is an isomorphism.

\goodbreak
\B{prop}\label{M-V-H} {\rm (Mayer--Vietoris for mixed
Hodge structures)\,} Associated with any square
(\ref{birsquare}) we have a long exact sequence of
mixed Hodge structures
$$\cdots\to H^n(X,\Z)\by{u}H^n(X',\Z)\oplus H^n(Y,\Z)\by{v}
H^n(Y',\Z)\by{\partial} H^{n+1}(X,\Z)\to\cdots $$ in which
$$ u=\left(\begin{array}{c}{f^*}\\{i^*}\end{array} \right)
\qquad{\rm and}\qquad v= (i'^*,-f'^*)$$
\E{prop}
\B{proof} We have a map of long exact sequences
$$ \begin{array}{ccccccccc}
\cdots & \to & H^n(X{\rm mod\,} Y,\Z) &
\to & H^n(X,\Z) & \by{i^*} &
H^{n}(Y,\Z) & \to & \cdots \\
 & &\p{f^*}\downarrow\cong &
&{\p{f^*}\downarrow}  & & {\p{f'^*}\downarrow} & &\\
\cdots & \to & H^n(X'{\rm mod\,} Y',\Z) &
\to & H^n(X',\Z) & \by{i'^*} &
H^{n}(Y',\Z) & \to & \cdots
\end{array} $$
where $H^*(-{\rm mod\,}\dag ,\Z)$ is the relative singular
cohomology functor (defined in \cite[8.3.8]{D}).
By excision $f^*\colon H^*(X{\rm \,mod\,}Y,\Z)\cong
H^*(X'{\rm\,mod\,}Y',\Z)$ (\cf \cite[8.3.10]{D}).
By \cite[8.3.9 and 8.2.2]{D} the diagram above is a diagram
in the abelian category of mixed Hodge structures.
The Mayer--Vietoris  exact sequence now follows by
a standard diagram chase.
\E{proof}
We then have, as well:

\B{schol} \label{M-V-D}
{\rm (Mayer--Vietoris for Deligne cohomology)\,}
Associated with any square (\ref{birsquare})
we have a long exact sequence in Deligne cohomology
$$\cdots\to H_{\cD}^n(X',\Z(r))\oplus
H_{\cD}^n(Y,\Z(r)) \to H_{\cD}^n(Y',\Z(r))\to
H_{\cD}^{n+1}(X,\Z(r))\to\cdots $$
\E{schol}
\B{proof} The proof of Proposition~\ref{M-V-H} goes through,
once we know that Deligne cohomology satisfies excision.
But since we have excision for the mixed Hodge
structure on relative singular cohomology, one can see it
holds for Deligne cohomology by arguing with the
relative cohomology sequence (\ref{modf-rel}).
\E{proof}

\B{thm}\label{M-V-K}
{\rm (Mayer--Vietoris for $K$-theory)\,}
Let $X$ be a reduced quasiprojective surface over a field
with normalization $\tilde X$.
Then there is a 1-dimensional subscheme $Y$ with
$Y_{red}={\rm Sing}\, X$ such that the
normalization square (\cf (\ref{birsquare}))
$$\B{array}{rcl} \tilde Y& \into &
\tilde X\\ \p{\tilde \pi}\downarrow & &\downarrow\p{\pi}\\
Y&\into & X \E{array}$$
induces exact sequences in $K$-theory:
$$
K_1(\tilde X)\oplus K_1(Y)\to K_1(\tilde Y)\by{\partial}
K_0(X)\to K_0(\tilde X)\oplus K_0(Y) \to K_0(\tilde Y)
$$
$$
SK_1(\tilde X)\oplus SK_1(Y)\to SK_1(\tilde Y)\by{\partial}
SK_0(X)\to SK_0(\tilde X)\to 0
$$
\E{thm}
\B{proof}
Let $K_*(X,\tilde X)$ and $K_*(Y,\tilde Y)$ be the
relative groups fitting into the long exact sequences in the
commutative diagram
$$\begin{array}{ccccccccccc}
K_1(X) &\to& K_1(\tilde X) &\to& K_0(X,\tilde X) &\to& K_0(X) &\to&
K_0(\tilde X)  &\to& K_{-1}(X,\tilde X)\\
\downarrow && \downarrow&&\downarrow &&\downarrow && \downarrow&&\veq\\
K_1(Y) &\to& K_1(\tilde Y) &\to& K_0(Y,\tilde Y) &\to& K_0(Y) &\to&
K_0(\tilde Y)  &\to& K_{-1}(Y,\tilde Y). \\
\end{array}$$
(The far right terms are isomorphic by \cite[A.6]{PW3}.)
To establish the existence and exactness of the $K_1-K_0$ sequence
we must show that ``excision'' holds for $K_0$, \ie that
$K_0(X,\tilde X)\cong K_0(Y,\tilde Y)$ for some $Y$
with $Y_{red}={\rm Sing}\, X$  (see \cite[5.1]{GW}).
If $Y$ is a subscheme of $X$ defined by an $\cO_{\tilde X}$-ideal
$\cI\subset\cO_X$
then by \cite[A.6]{PW3} there is a natural exact sequence
\B{equation}\label{K0rel}
H^1(Y,{\cI}/{\cI}^2\otimes\Omega_{\tilde X/X})
\vlongby{\eta(Y)}
K_0(X,\tilde X) \to K_0(Y,\tilde Y) \to0.
\E{equation}
We define $Y_1$ using the conductor ideal $\cJ$, and $Y$ using
the ideal ${\cI}={\cJ}^2$.
Then $Y_{red}={\rm Sing}\, X$, and the map from
${\cal I}/{\cal I}^2$ to ${\cal J}/{\cal J}^2$ is zero.
By naturality in $Y\to Y_1$, the map $\eta(Y)$ in (\ref{K0rel})
is the composite map
$$
  H^1(Y,{\cal I}/{\cal I}^2\otimes\Omega_{\tilde X/X}) \by0
  H^1(Y,{\cal J}/{\cal J}^2\otimes\Omega_{\tilde X/X})
\vlongby{\eta(Y_1)}  K_0(X,\tilde X)$$
so $\eta(Y)=0$ in (\ref{K0rel}).
Hence excision holds for $Y$, as claimed.

There is a natural map from the $K_1-K_0$ sequence onto the
 ``Units-Pic'' sequence, and the kernel is the $SK_1-SK_0$ sequence.
A standard diagram chase, described in  \cite[8.6]{PW1},
shows that the latter sequence sequence is also exact.
\E{proof}

We remark that if $Y$ is reduced, or 0-dimensional, or even affine, then
the obstruction $H^1(Y)$ in (\ref{K0rel}) automatically vanishes,
and excision is immediate.  Theorem \ref{M-V-K} was proven in these
special cases in \cite[7.5]{PW1} and \cite[A.3]{PW3}.

\B{cor}\label{KtoD}
With the notation of \ref{M-V-K}, the following diagram commutes.
$$
\begin{array}{cccccccc}
SK_1(\tilde X)\oplus SK_1(Y) &\kern-8pt\to\kern-9pt&
SK_1(\tilde Y) &\kern-7pt\to\kern-8pt &SK_0(X) &
\kern-7pt\to\kern-8pt & SK_0(\tilde X) &\kern-7pt\to0\\
\downarrow&&\downarrow&&\downarrow&&\downarrow&\\
K_1(\tilde X)\oplus K_1(Y) &\kern-8pt\to\kern-9pt&
K_1(\tilde Y) &\kern-7pt\to\kern-8pt &K_0(X) &
\kern-7pt\to\kern-8pt & K_0(\tilde X) &\kern-7pt\to0\\
c_2\downarrow\ &&c_2\downarrow\ &&c_2\downarrow\ &&c_2\downarrow\ &\\
H_\cD^3(\tilde X,\Z(2))\oplus H_\cD^3(Y,\Z(2)) &
\kern-8pt\to\kern-9pt & H_\cD^3(\tilde Y,\Z(2))&
\kern-7pt\to\kern-8pt &
H_\cD^4(X,\Z(2)) &\kern-7pt\to\kern-8pt &
H_\cD^4(\tilde X,\Z(2)) &\kern-7pt\to0
\end{array}
$$
\E{cor}
\B{proof} We use the notation of (2.4.1). 
For each open $U$ in $X$,
let $F(U)$ denote the homotopy fiber of
$K(U\times_X Y)\times K(U\times_X\tilde X)
\to K(U\times_X\tilde Y)$.  By Proposition \ref{M-V-D} the
corresponding homotopy fiber for Deligne cohomology is
$\cD(U)$.  In addition, there is a natural map from
$K(U)$ to $F(U)$ which is an isomorphism on $\pi_0$ by
Theorem \ref{M-V-K}.  Therefore the natural map $C_2^{ss}$
of (2.4.1) 
induces a map
$F(U)\to\cD(U)$ on homotopy fibers, making the diagram
$$\begin{array}{ccccc}
K_1(\tilde Y) &\by{\partial}& \pi_0F(X)
&{\stackrel{\simeq}{\leftarrow}}& K_0(X)\\  \strut
C_2^{ss}\downarrow\ &&\ \downarrow C_2^{ss}&&\ \downarrow C_2^{ss}\\
H_\cD^3(\tilde Y,\Z(2)))&\by{\partial}& \pi_0\cD(X)
&{\stackrel{\simeq}{\leftarrow}}& H_\cD^4(X,\Z(2))
\end{array}$$
commute. But the top composite is the $K$-theory boundary map
in Theorem \ref{M-V-K}.
\E{proof}
Using (2.4) and Lemma~\ref{c2-SK1}, we may refine
Corollary~\ref{KtoD} as follows.
\B{schol}\label{SKtoD}
With the notation of \ref{M-V-K}, the following diagram commutes.
$$
\begin{array}{cccccccc}
SK_1(\tilde X)\oplus SK_1(Y) &\kern-8pt\to\kern-9pt&
SK_1(\tilde Y) &\kern-7pt\to\kern-8pt &SK_0(X) &
\kern-7pt\to\kern-8pt & SK_0(\tilde X) &\kern-7pt\to0\\
\mbox{\rm onto}\downarrow\ \ &&\cong\downarrow&&\cong\downarrow&&
\cong\downarrow&\\
H^1(\tilde X,\cK_2)\oplus H^1(Y,\cK_2) &\kern-8pt\to\kern-9pt&
H^{1}(\tilde Y,\cK_2) &\kern-7pt\to\kern-8pt &
H^{2}(X,\cK_2) & \kern-7pt\to\kern-8pt &
H^{2}(\tilde X,\cK_2)& \\
c_2\downarrow\ &&c_2\downarrow\ &&c_2\downarrow\ &&c_2\downarrow\ &\\
H_\cD^3(\tilde X,\Z(2))\oplus H_\cD^3(Y,\Z(2)) &
\kern-8pt\to\kern-9pt & H_\cD^3(\tilde Y,\Z(2))&
\kern-7pt\to\kern-8pt &
H_\cD^4(X,\Z(2)) &\kern-7pt\to\kern-8pt &
H_\cD^4(\tilde X,\Z(2)) &\kern-7pt\to0
\end{array}
$$
\E{schol}

\goodbreak
\section{The Albanese 1-motive of a proper surface}

In this section a {\it surface}\, will mean a
\underline{proper}
reduced $2$-dimensional scheme $X$ of finite type
over the complex numbers $\C$.
We will consider the intermediate jacobian $$J^2(X) \df
\frac{H^{3}(X,\C)}{F^2H^{3}+ H^{3}(X,\Z(2))}$$
This is the mixed Hodge theoretic
generalization of the classical Albanese group variety of a
smooth surface.

We begin with an elementary result
(cf. \cite[Remark 5.5]{GDuke}).

\B{lemma}\label{filt}
Suppose that $X$ is a proper surface.  Then
$$
F^2H^i(X,\C)\cap H^i(X,\R)=0 \mbox{\ \  for \ } i=2,3.
$$
Hence  in sequence (\ref{modf}) we have
$$\B{array}{ccc} H^i(X,\Z)_{tors}
\strut&=&\mbox{kernel of } H^i(X,\Z)\to H^i(X,\C)/F^2H^i\\
\strut&=&\mbox{image of } H_\cD^i(X,\Z(2))\by{\varepsilon}H_{an}^i(X,\Z).
\E{array}$$
\E{lemma}

\kern-6pt
\B{proof}
We will show  that $H^i(X,\R)$ injects into $H^i(X,\C)/F^2$.
When $X$ is smooth, then $H^i(X)$ has pure weight $i$.
In this case complex conjugation on $H^i(X,\C)$ fixes
$H^i(X,\R)$ but the subspace $F^2H^i(X,\C)$ meets its
conjugate in $0$.

If $X$ is a singular surface, choose a resolution of
singularities $X'\to X$.  If $Y$ is a curve containing
the singular locus of $X$, then we are in the situation
of square (\ref{birsquare}).  Since
$F^2H^1=F^2H^2=0$ for the curves $Y$ and $Y'$,
the Mayer--Vietoris sequence in Proposition~\ref{M-V-H}
yields $F^2H^i(X,\C)=F^2H^i(X',\C)$ for $i=2,3$.
Comparing the $\R$ and $\C$ structures in the
Mayer--Vietoris long  exact sequence of Proposition~\ref{M-V-H}
yields the following diagram has exact rows:
$$\begin{array}{ccccc}
H^1(Y',\R)&\to&H^2(X,\R) & \to & H^2(X',\R)\\
\downarrow & &\downarrow & &\downarrow \\
H^1(Y',\C)&\to& \displaystyle\frac{H^2(X,\C)}{F^2H^2}&
\to &\displaystyle\frac{H^2(X',\C)}{F^2H^2}
\end{array}$$
The right-most vertical arrow in the diagram
is injective because $X'$ is smooth.
A diagram chase shows that the middle vertical
arrow is injective, whence the lemma.
\E{proof}

\subsection*{Normal surfaces}

Consider a surface $X$ with normal singularities;
its singular locus $\Sigma$ is a finite set of closed points.
Choose a desingularization $f\colon X' \to X$ and consider the
exceptional divisor $E=f^{-1}(\Sigma)$;  $E$ is the finite disjoint
union of the inverse images $E_{\sigma}\df f^{-1}(\sigma)$
of the $\sigma\in \Sigma$.
Associated to $f$ is the square (\ref{birsquare}),
with $Y=\Sigma$ and $Y'=E$.  Because the fibers $E_\sigma$ of $f$
are connected (by Zariski's Main Theorem), we have
$H^0(\Sigma,\Z)\cong H^0(E,\Z)$.
 From Proposition~\ref{M-V-H} we get a
long exact sequence of mixed Hodge structures:
\B{equation}\label{exnorm}
\B{array}{c}
0\to H^1(X,\Z)\by{f^*} H^1(X',\Z)\to
H^1(E,\Z)\to H^2(X,\Z)\by{f^*} \qquad\\ \qquad\qquad
H^2(X',\Z)\to H^2(E,\Z)\to H^3(X,\Z)\by{f^*} H^3(X',\Z) \to 0
\E{array} \E{equation}
\noindent
Recall that if $X$ is proper then
each $H^n(X')$ has pure weight $n$, and
that $W_{n-1}H^n(X,\Q)$ is the kernel of
$f^*\colon H^n(X,\Q)\to H^n(X',\Q)$ by \cite[8.2.5]{D}.
That is,  $H^n(X)$ has pure weight $n$ if and only if
$H^n(X,\Q)$ injects into $H^n(X',\Q)$.
There are examples of normal surfaces for which $H^2(X)$
does not have pure weight 2, i.e.,  with $W_1H^2\neq 0$
(\cf \cite{BVS1},\cite{BVS2}).  The following result
quantifies this impurity.

\B{prop}\label{pure} Let $X$ be a proper normal surface.
If $n\ne2$ then $H^n(X)$ has pure weight $n$.  If $n=2$ and
$E$ is the exceptional divisor in a desingularization $X'$, then
$$W_{1}H^2(X,\Q)=\coker H^1(X',\Q)\to H^1(E,\Q).$$
\E{prop}
\B{proof} If $n\ne3$ this follows from the
sequence of mixed Hodge structures (\ref{exnorm}).
For $n=3$ we must show that $H^3(X,\Q)$ embeds in
$H^3(X',\Q)$.  Nothing is lost if we replace $X'$
by a quadratic transformation, so
we may assume that all the irreducible components
$E_1,\dots,E_n$ of $E$ are non-singular and that
if $i\neq j$ and $E_i\cap E_j \neq\emptyset$ then
$E_i$ and $E_j$ intersect transversally in exactly one point
not belonging to any other $E_k$.
 From (\ref{exnorm}) and the commutative square
$$\B{array}{ccc} \Pic(X') &\by{\cap E} & \Pic(E) \\
\p{c_1}\downarrow& & \downarrow \p{c_1}\\
H^2(X',\Z) &\to & H^2(E,\Z)
\E{array}$$
we see that it suffices to prove that
$\Pic(X')\tensor\Q \to H^2(E,\Q)$ is a surjection.
Now $H^2(E,\Q)=\oplus H^2(E_i,\Q) \cong \Q^n$,
and  $\Pic(E)\otimes\Q  \by{c_1} H^2(E,\Q)$
is just the degree map.
Moreover, the intersection pairing on the divisors
on $X'$ satisfies $(D\bul E_i) ={\rm deg}(D\cap E_i)$.
Thus if we represent an element of
$\Pic(X')$ by a divisor $D$, its image in
$H^2(E,\Q)\cong\Q^n$ is given by the intersection vector
$(D\bul E_1,\dots,D\bul E_n)$.
Now each $E_i$ represents an element of $\Pic(X')$,
and their intersection vectors form a basis of $H^2(E,\Q)$
because the intersection matrix
$(E_i\bul E_j)$ is negative definite
(see \cite[\S 1]{M} or \cite[14.1]{Lip}).
\E{proof}

\B{cor}\label{abnorm} Let $f\colon X' \to X$ be a
resolution of singularities of a proper normal surface. Then
$J^2(X)$ is an abelian variety, because there is an isomorphism
$$f^*\colon J^2(X)\by{\simeq} J^2(X')$$
\E{cor}
\B{proof} By (\ref{exnorm}) and \ref{pure},
$f^*\colon H^3(X,\Z)\by{} H^3(X',\Z)$
is onto with torsion kernel.  \E{proof}

\subsection*{Normalization}

Now let $X$ be a non-normal surface.
The singular locus $\Sigma$ of $X$ is
$1$-dimensional. 
Letting $\pi: \tilde X \to X$ denote its normalization,
we have $\pi\colon\tilde X-\tilde\Sigma\by{\simeq}X-\Sigma$,
where $\tilde\Sigma=\pi^{-1}\Sigma$.
By Proposition~\ref{M-V-H},
$\pi$ induces a long exact 
sequence of mixed Hodge structures.
\B{equation}\label{exnon-norm}
\B{array}{c}
H^1(X,\Z)\to H^1(\tilde X,\Z)\oplus H^1(\Sigma
,\Z)\to  H^1(\tilde\Sigma ,\Z)\to H^2(X,\Z)\to\\
H^2(\tilde X,\Z)\oplus H^2(\Sigma,\Z)\to H^2(\tilde\Sigma ,\Z)
\to H^3(X,\Z) \to H^3(\tilde X,\Z) \to 0  \E{array}
\E{equation}
Since the Hodge structure on $H^2$ of a curve is pure
of type $(1,1)$, the abelian group
$$M = \frac{\coker H^2(\Sigma,\Z)\by{\pi^*}H^2(\tilde\Sigma,\Z)}
           {\coker H^2(X,\Z)\by{\pi^*}H^2(\tilde X,\Z)}$$
has a mixed Hodge structure which is pure of type $(1,1)$,
and there is an extension of mixed Hodge structures
\B{equation}\label{exMHS}
0\to M \to H^3(X,\Z) \to H^3(\tilde X,\Z) \to0.
\E{equation}

\B{prop}\label{non-norm} Let $X$ be a proper surface,
with normalization $\pi\colon\tilde X \to X$.
Then we have an extension
$$0\to (\C/\Z)^s \to J^2(X)\by{\pi^*} J^2(\tilde X)\to 0$$
where $s$ is the rank of the abelian group $M$.
\E{prop}
\B{proof} $J^2(X)$ is the cokernel of the natural map
$H^3(X,\Z(2))\to H^3(X,\C)/F^2H^3$.  Given this, the result
is a formal consequence of (\ref{exMHS})
and the fact that $F^2M=0$, which implies
that $F^2H^3(X,\C) \cong F^2H^3(\tilde X,\C)$.
We remark that the complex torus $(\C/\Z)^s$ that arises in this
extension is a quotient of the complex torus
$(\C/\Z(2))^s = M\otimes(\C/\Z(2))$ by a finite group.
\E{proof}
\B{cor}\label{1-motive}
Let $f\colon X'\to X$ be a desingularization of a proper
surface $X$, obtained by resolving the singularities of its
normalization $\tilde X$.
Then there is an exact sequence
$$0\to (\C/\Z)^s \to J^2(X)\by{f^*} J^2(X')\to 0$$
where $s$ is the rank of $M$, as in Proposition~\ref{non-norm}.
In particular, if $X'$ has irregularity $q$ then the torsion
subgroup of $J^2(X)$ is isomorphic to $(\Q/\Z)^{2q+s}$.
\E{cor}

Recall from  \cite[10.1.2]{D} that a ``1-motive''
 $M=(L,A,T,J,u)$ is defined to be an extension $J$
of an abelian variety $A$ by a complex torus $T$, a lattice
$L$ and a homomorphism $L\by{u} J$.  Since we may
canonically identify the group of $\C$-points of the
abelian variety $Alb(X')$ with $J^2(X')$,
the conclusion of Corollary~\ref{1-motive}
is just that $J^2(X)$ is part
of a 1-motive $Alb(X)$ in which the lattice $L$ is zero.


\B{defi} Let $X$ be a proper surface over $\C$.
The {\it Albanese 1-motive} of $X$ is the 1-motive $Alb(X)$
given by $$(0,Alb(X'),(\C/\Z(2))^s,J^2(X),{\rm zero}).$$
As the construction in \ref{non-norm} shows,
$Alb$ is a functor from proper surfaces to 
1-motives.
\E{defi}

\goodbreak
\subsection*{Torsion in $J^2(X)$}

For simplicity, let us write $\Q/\Z$ for the torsion subgroup
$\Q(2)/\Z(2)$ of $\C/\Z(2)$, so that
$H^i(-,\Q/\Z)\cong H^i(-,\C/\Z(2))_{tors}$.
The maps $H^i(-,\C/\Z(2))\to H_\cD^{i+1}(-,\Z(2))$
of (\ref{star}) induce canonical maps
$$H^i(-,\Q/\Z)\cong
H^i(-,\C/\Z(2))_{tors}\to H_\cD^{i+1}(-,\Z(2))_{tors}.$$
These are the maps in the following Proposition:

\B{prop}\label{HD-tors}
Let $Z$ be a proper scheme over $\C$.
Then
\B{description} 
\item[{\it i)}] $H^1(Z,\Q/\Z)\by{\cong} H_\cD^2(Z,\Z(2))_{tors}$
\item[{\it ii)}] If $Z$ is either a curve or a surface then
$$H^2(Z,\Q/\Z)\by{\cong} H_\cD^3(Z,\Z(2))_{tors} $$
%
\item[{\it iii)}] If $Z$ is a surface then
$$H^3(Z,\Q/\Z)\by{\cong}H_\cD^4(Z,\Z(2))_{tors}\cong
 J^2(Z)_{tors}$$
\E{description} 
\E{prop}
\B{proof} The first assertion was proven in \ref{sheafHD}.
If $Z$ is a curve then $F^2H^2(Z,\C)=0$, so by
(\ref{star}) we have $H_\cD^3(Z,\Z(2))\cong H^2(Z,\C/\Z(2))$,
and the result is immediate.

We may therefore suppose that $Z$ is a surface, say
with $c$ irreducible components, so that $H^4(Z,\Z)=\Z^c$.
We deduce from (\ref{extJZ})  that
$J^2(Z)_{tors}\cong H_\cD^4(Z,\Z(2))_{tors}$.
Moreover, since $F^2H^4(Z,\C)=H^4(Z,\C)=\C^c$
the sequence (\ref{star}) ends in
\B{equation}\label{star3}  H_\cD^3(Z,\Z(2))\to
F^2H^3(Z,\C)\to H^3(Z,\C/\Z(2))\to H_\cD^4(Z,\Z(2))
\by{\varepsilon}\Z^c\to0.  \E{equation}
Lemma \ref{filt} states that for $i=2,3$ the image of
$\varepsilon$ in the exact sequence
$$H^{i-1}(Z,\C)/F^2H^{i-1}\to H_\cD^i(Z,\Z(2))\by{\varepsilon}
H^i(Z,\Z)\to H^i(Z,\C)/F^2H^i$$
of (\ref{modf}) is the torsion subgroup $H^i(Z,\Z)_{tors}$.
Combining this with  the universal coefficient theorem,
we have a commutative diagram with exact rows:
$$\begin{array}{ccccccc}
0\to\kern-7pt&H^{2}(Z,\C)/F^2H^{2}\oplus H^{2}(Z,\Z)
&\to&H_\cD^3(Z,\Z(2)) &\by{\varepsilon}&H^3(Z,\Z)_{tors}&\to0\\
 &\uparrow & &\uparrow & &\veq & \\
0\to\kern-7pt&H^{2}(Z,\Z)\otimes\C/\Z(2) &\to &
H^{2}(Z,\C/\Z(2))&\to & H^3(Z,\Z)_{tors}&\to 0.
\end{array}$$
By the five-lemma, sequence (\ref{star}) and (\ref{star3})
we get the extensions
\B{equation}
0\to F^2H^2(Z,\C)\to H^{2}(Z,\C/\Z(2))\to H_\cD^3(Z,\Z(2))\to0,
\E{equation}
\B{equation}
0\to F^2H^3(Z,\C)\to H^3(Z,\C/\Z(2))\to J^2(Z)\to0.
\E{equation}
Since $F^2H^i(Z)$ is uniquely divisible, we may
pass to torsion subgroups.
This proves the remainder of the proposition.
\E{proof}
\medskip



\section{Curves} 

The singular locus of a reduced surface is usually an
(unreduced) curve. For this reason, we need information
about $K_1$ and $K_2$ of curves in order to study
surfaces. This information is given by theorems
\ref{smooth} and \ref{sing} below.
Part {\it i)\,} of Theorem \ref{smooth}
is of course well-known and
almost classical; a reference is \cite[1.1]{R}.
Since these results are of independent
interest, we have expanded our exposition to include
the case of characteristic $p$.

By a `curve' over a field $k$ we mean a $1$-dimensional
quasiprojective scheme $Y$ over $k$;
a curve is not necessarily reduced.
There is a natural map from $K_1(Y)$ to the
group  $H^0(Y,\cO_Y)$ of global units of $Y$; the kernel
of this map is usually written as $SK_1(Y)$.
When $Y$ is a curve there is a natural isomorphism
$SK_1(Y)\cong H^1(Y,\cK_2)$, as well as
a natural short exact sequence
$$0\to H^1(Y,\cK_3)\to K_2(Y)\to H^0(Y,\cK_2)\to 0$$
given by the Brown-Gersten spectral sequence \cite{TT}.

\B{thm}\label{smooth} Let $Y$ be a smooth curve over an
algebraically closed field $k$. Let $r\geq 0$ denote the
number of irreducible components of $Y$ which are proper.
Then \B{description}
\item[{\it i)}] $SK_1(Y)\cong(k^*)^r\oplus V_1$
where $V_1$ is a uniquely divisible group;
\item[{\it ii)}] $K_2(Y)$ and $H^0(Y,\cK_2)$ are both
divisible abelian groups.
\E{description}
 \E{thm}

\noindent\B{proof}
If $Y$ is a smooth connected curve over an algebraically
closed field $k$ then the localization sequence is
\B{equation}\label{localization}
\coprod_{y\in Y(k)}^{} K_2(k) \to K_2(Y) \to
K_2(k(Y))\by{tame} \coprod_{y\in Y(k)}^{} k^* \to
SK_1(Y)\to 0
\E{equation}
and the image of $K_2(Y)$ in $K_2(k(Y))$ is
$H^0(Y,\cK_2)$. Since $\coprod K_2(k)$
is divisible \cite[1.3]{BT}, the divisibility of $K_2(Y)$
is equivalent to the divisibility of $H^0(Y,\cK_2)$. If
char$(k)=0$, the result now follows from Suslin's exact
sequence \cite[4.4]{S2} for $n$ invertible in $k$:
$$0\to H^0(Y,\cK_2)/n\to \Het{2}(Y,\mun)\to {}_nSK_1(Y) \to0$$
Indeed, if $Y$ is affine then $\Het{2}(Y)=0$,
and if $Y$ is projective then the composite $\mu_n \cong
\Het{2}(Y,\mun)\to SK_1(Y)\to k^*$
is the standard inclusion. 

If char$(k)=p>0$, we need only a slight additional argument.
Because $k(Y)$ is the function field of a curve,
we know from \cite[p.391]{BT} that $K_2(k(Y))$ is
$p$-divisible, and from \cite[1.10]{S2}
(which is implicit in \cite[p.397]{BT}) that it
has no $p$-torsion. Hence both $K_2(k(Y))$ and
$\coprod k^*$ are uniquely $p$-divisible groups, \ie
$\Z[\frac{1}{p}]$-modules. It follows that both the
kernel $H^0(Y,\cK_2)$ and cokernel $SK_1(Y)$ of the
`tame symbol' map in (\ref{localization})  must be
uniquely $p$-divisible. This proves Theorem~\ref{smooth}
in characteristic $p$.
\E{proof}

\B{lemma} Let $Y$ be a smooth connected projective curve
over $\C$.  Then
$$c_2\colon SK_1(Y) \to H_\cD^3(Y,\Z(2))\cong \C^*$$
is a split surjection.
In particular, it is an isomorphism on torsion subgroups.
\E{lemma}
\B{proof} This is implicit in p.219 of Gillet's paper \cite{GDuke}.
The isomorphism $H_\cD^3(Y,\Z(2))\cong \C/\Z(2)\cong\C^*$
follows from (\ref{star}) or \ref{HD-tors}.
If $y\in Y(\C)$ is considered as an
element of $Pic(Y)$ and $z\in\C^*$ then we can form
$\{ y,z\}\in SK_1(Y)$ and the product formula yields
$c_2(\{ y,z\})=-c_1(y)\cup c_1(z) = z^{-1}\in \C^*$.
\E{proof}

\B{thm}\label{sing} Let $Y$ be any curve over an
algebraically closed field $k$, 
and let $r\geq 0$ denote the number of irreducible
components of $Y$ which are proper. Then \B{description}
\item[{\it i)}] If char$(k)=0$, or if $Y$ is reduced, then
$$SK_1(Y)\cong (k^*)^r\oplus V_1,$$
where $V_1$ is a uniquely divisible abelian group;
\item[{\it ii)}] If char$(k)=p>0$ then
$$SK_1(Y)\cong (k^*)^r\oplus V_1\oplus P,$$
where $V_1$ is uniquely divisible and $P$ is a
$p$-group of bounded exponent.
\item[{\it iii)}] If $k=\C$  then the Chern class
$$c_2\colon SK_1(Y) \to H_\cD^3(Y,\Z(2))\cong (\C^*)^r$$
is a split surjection.
In particular, it is an isomorphism on torsion subgroups.
\E{description}
 \E{thm}

\noindent\B{proof}
We proceed in three steps.\\
{\it Step 1.}
Suppose that $Y$ is any reduced curve over $k$. If we
pick $r$ smooth points $y_i$ on $Y$, one on each
proper component of $Y$, then $Y_0=Y-\{y_1,\ldots ,y_r\}$
is affine. The localization sequence for $Y_0\subset Y$ is
\B{equation}\label{notloc}
\coprod_{i=1}^{r} K_2(k) \to K_2(Y) \to
K_2(Y_0)\to \coprod_{i=1}^{r} k^* \to
SK_1(Y)\to SK_1(Y_0)\to 0
\E{equation}
If $\tilde Y$ is the normalization of $Y$, then we may
indentify the $y_i$ with points on the smooth curve
$\tilde Y$. By the smooth case \ref{smooth}, the composition of
$$\coprod k^* \to SK_1(Y)\to SK_1(\tilde Y)$$
is an injection. Hence $SK_1(Y)$ is the direct sum of
$\coprod k^*$ and $SK_1(Y_0)$ while $K_2(Y)$ is the direct
sum of the image of the divisible group $\coprod K_2(k)$
and the group $K_2(Y_0)$.
Part (iii) now follows from the above Lemma.
This argument also shows that we may replace $Y$ by
$Y_0$ in proving parts (i) and (ii) of Theorem~\ref{sing}
for reduced curves.\\
{\it Step 2.} Now suppose that $Y= \mbox{Spec}(A)$ is any
reduced affine curve over $k$. Let $B$ be the
normalization of $A$, and $I$ the conductor ideal from
$B$ to $A$. By \cite[3.1 and 4.2]{GR}, excision holds for
$K_1$ and there is an exact sequence
$$K_2(B)\oplus K_2(A/I)\to K_2(B/I)\to SK_1(A)\to
SK_1(B)\to 0$$
Since $B$ is a finite product of Dedekind domains, $B/I$
is a finite principal ideal ring. By Corollary~\ref{k2div}
below, $K_2(B/I)$ is uniquely divisible.
By Theorem~\ref{smooth}, $SK_1(B)$ is uniquely divisible
and $K_2(B)$ is divisible. Finally, since  $A/I$ is finite
dimensional, we know from Corollary~\ref{k2div} that
$K_2(A/I)$ is uniquely divisible (modulo bounded $p$-torsion if
char$(k)=p\neq 0$). A diagram chase shows that $SK_1(A)$ is
uniquely divisible (modulo bounded $p$-torsion if $p\ne0$).
This proves theorem \ref{sing} for reduced curves.

\B{lemma}\label{nilpo} Let $I$ be a nilpotent ideal in an
algebra $A$ over a field $k$.
\B{description} \item[{\it (a)}] If char$(k)=0$,
$K_n(A,I)$ is a uniquely divisible group for every $n$.
\item[{\it (b)}] If char$(k)=p>0$, $K_2(A,I)$ is a
$p$-group of bounded exponent.
 \E{description}
\E{lemma}
\B{proof} Part {\it (a)\,} is proven in
\cite[1.4]{W966}.  
If char$(k)=p$, chose $m$ such that
$I^{p^m}=0$; we will show that $p^mK_2(A,I)$. Indeed,
$K_2(A,I)$ is generated by Steinberg symbols $\{a,1+x\}$
with $a\in A$ and $x\in I$, and $p^m\,\{a,1+x\}$ is
$\{a,1+x^{p^m}\}=\{a,1\}=0$.
 \E{proof}
 \B{cor}\label{k2div} Let $A$ be a finite algebra over
an algebraically closed field $k$.
\B{description} \item[{\it (a)}] If char$(k)=0$ or if $A$
is a principal ideal ring, the group $K_2(A)$ is uniquely divisible.
\item[{\it (b)}] If char$(k)=p$, $K_2(A)$  is the sum of
the uniquely divisible group $K_2(A_{red})$ and a
$p$-group of bounded exponent.
 \E{description}
 \E{cor}
\B{proof} Let $I$ be the nilradical of $A$, so that
$A_{red}=A/I$ is semisimple and hence $A\to A_{red}$
splits. Then $K_2(A)\cong K_2(A_{red})\oplus K_2(A,I)$,
and $K_2(A_{red})$ is uniquely divisible by
\cite[1.3]{BT}. Finally, if $A$ is a principal ideal
ring then $A$ is a product of truncated polynomial rings
$k[s]/(s^n)$ and a direct calculation (\cite[p.485]{Graham})
shows that $K_2(k[s]/(s^n))\cong K_2(k)$.
 \E{proof}

\noindent{\it Step 3.} Finally, suppose that $Y$ is a curve
which is not reduced.
Let $\cI$ denote the nilradical
ideal sheaf of $\cO_Y$, and write $\cK_2\cI$ for the
sheafification of the presheaf
$U\mapsto K_2(\cO_Y(U),\cI(U))$.
If $\cK_{2,red}$ denotes the sheafification of
$U\mapsto K_2(U_{red})$, there is an
exact sequence of sheaves on $Y_{Zar}$
\B{equation}\label{Ired}
\cK_2\cI\to \cK_2\to \cK_{2,red}\to 0
\E{equation}
Let $U$ denote the smooth locus of $Y_{red}$. Since
$U_{red}$ is smooth, the ring map $\cO_U\to\cO_{U_{red}}$
splits. Therefore $\cK_2\cI\mid_U$ injects into
$\cK_2\mid_U$, \ie the kernel of $\cK_2\cI\to\cK_2$ is a
skyscraper sheaf supported on $Y-U$. It follows that we
have an exact sequence
$$H^0(Y,\cK_{2,red})\to H^1(Y,\cK_2\cI)\to H^1(Y,\cK_2)\to
H^1(Y_{red},\cK_2)\to 0$$
which we may rewrite as follows
\B{equation}\label{reduced}
K_2(Y_{red})\to H^1(Y,\cK_2\cI)\to SK_1(Y)\to
SK_1(Y_{red})\to 0
\E{equation}
By Step 2, $SK_1(Y_{red})$ is uniquely divisible. If
char$(k)=p$, we know by Lemma~\ref{nilpo}(b) that
$H^1(Y,\cK_2\cI)$ is a $p$-group of bounded exponent,
and this part {\it (ii)\,} of Theorem~\ref{sing} because
a uniquely divisible group has no nontrivial extensions
by a $p$-group.  Finally, suppose that char$(k)=0$.
By Lemma~\ref{nilpo}(a), $H^1(Y,\cK_2\cI)$ is uniquely
divisible. By Proposition~\ref{n-divisible} below,
$K_2(Y_{red})$ is divisible.
In this case part {\it (i)\,} of Theorem~\ref{sing}
follows from Step 2 and the exact sequence (\ref{reduced}).
\E{proof}

\B{prop}\label{n-divisible} If $\frac{1}{n}\in k$ and
$Y$ is a curve then $K_2(Y)$ is $n$-divisible.
\E{prop}
\B{proof} We consider the $K$-theory of $Y$ with
coefficients $\Z/n$, which is related to the usual
Quillen $K$-theory of $Y$ by exact sequences such
as
$$0\to K_2(Y)\otimes\Z/n\to K_2(Y;\Z/n)\to
K_1(Y)_{n-tors}\to 0$$
We know by \cite[1.4]{W966} that $K_2(Y;\Z/n)\cong
K_2(Y_{red};\Z/n)$, and hence that $K_2(Y)\otimes\Z/n$ is
a subgroup of $K_2(Y_{red})\otimes\Z/n$. Thus we may
assume that $Y$ is reduced. Let $\tilde Y$ be the
normalization of $Y$. The conductor ideal defines a
zero-dimensional subscheme Spec($C$) of $Y$, and also its
preimage Spec($D$) in $\tilde Y$. Because excision holds
(see \cite[1.2]{W966}) we have an exact sequence
$$K_3(D;\Z/n)\to K_2(Y;\Z/n)\to K_2(\tilde Y;\Z/n)\oplus
K_2(C;\Z/n)$$
Now $K_3(D;\Z/n)\cong K_3(D_{red};\Z/n)$, again by
\cite[1.4]{W966}. Since $D_{red}$ is a finite product of
copies of $k$, and $K_3(k;\Z/n)=0$ by Suslin \cite{Suslin},
we
have $K_3(D;\Z/n)=0$. Hence $K_2(Y;\Z/n)$ injects into
$K_2(\tilde Y;\Z/n)\oplus K_2(C;\Z/n)$. By naturality,
the subgroup $K_2(Y)\otimes\Z/n$ of $K_2(Y;\Z/n)$ injects
into the corresponding subgroup $K_2(\tilde Y)\otimes\Z/n
\oplus K_2(C)\otimes\Z/n$ of $K_2(\tilde Y;\Z/n)\oplus
K_2(C;\Z/n)$, but: $K_2(\tilde Y)$ is divisible by
Theorem~\ref{smooth} and $K_2(C)$ is divisible by
Corollary~\ref{k2div}, so this latter subgroup is zero,
hence $K_2(Y)\otimes\Z/n =0$ as claimed.
 \E{proof}

\section{$K$-theory results}

In this section we collect some results on the relation
between the Zariski sheaves $\cK_2$ and $\cH^q(\mun)$,
namely \ref{HoobH2}, \ref{nK2} and \ref{square},
which will be used in the proof of the main theorem.
In this section, our field $k$ will always contain $\frac1n$.

The first result, which we cite without proof,
concerns the sheafification of the \'etale Chern class
$c_2^{et}\colon K_2(X)/n\to\Het2(X,\mun)$.
It is an extension by Hoobler of a
well known result for smooth schemes.

\B{prop}\label{HoobH2}
(Hoobler \cite{Hoob}; \cf \cite[Thm. 0.2]{PW3})
Let $X$ be a scheme of finite type over a
field containing $\frac1n$.
Then the \'etale Chern class $c_2^{et}$ induces
an isomorphism of Zariski sheaves on $X$: \quad
$\cK_2/n \by{\simeq} \cH^2(\mun)$.
\E{prop}

Our other results concern the n-torsion subsheaf
${}_n\cK_2$  of $\cK_2$.  We begin with the local version.

\B{lemma}\label{HoobH1} Let $A$ be a semilocal ring
essentially of finite type over a field $k$.
Assume $k$ contains a primitive $n^{th}$ root of unity $\zeta$.
Define a map $$\varphi:
\Het{1}(A,\mun)\cong A^*/A^{*n}\to{}_nK_2(A)$$ by
$\varphi(a)=\{a,\zeta\}$, $a\in A^*$. Then
$\varphi$ is surjective. If $A$ is regular and $k$ contains
an algebraically closed field
then  $\varphi$ is an isomorphism.
\E{lemma}
\B{proof} $\varphi$ is well defined because
$\{ a^n,\zeta\}=\{ a,1\}=0$.  Suppose first that $A$ is regular.
Then $\varphi$ is onto by the Merkurev-Suslin Theorem.
If in addition its field of fractions $F$ contains an
algebraically closed field, then
$\Het{1}(F,\mun)\cong{}_nK_2(F)$
by \cite[3.7]{S2}.  Comparing the Bloch-Ogus resolution of
$\Het{1}(A,\mun)$ to the Gersten--Quillen
resolution of ${}_nK_2(A)$, one gets that
$\Het{1}(A,\mun)\cong {}_nK_2(A)$.

The promotion to any semilocal ring $A$ follows from the same
arguments used by R.Hoobler in \cite{Hoob}. Since $A$ is a
localization of a finitely generated $k$-algebra, there
exists localization $B$ of a polynomial ring over $k$
and an ideal $I$ in $B$ such that $A=B/I$.
Let $(B^h,I^h)$ be the henselization of the pair $(B,I)$.
As $B^h$ is a direct limit of semilocal regular rings
finite over $B$,
the map $\Het{1}(B^h,\mun)\by{\varphi}{}_n\cK_2(B^h)$
is an isomorphism. By a result of O. Gabber \cite[Th.1]{Gabber}
we have $K_3(B^h;\Z/n)\cong K_3(B^h/I^h;\Z/n)\cong K_3(A;\Z/n)$.
By proper base change
$$\Het{1}(B^h,\mun)\cong
\Het{1}(B^h/I^h,\mun)\cong
\Het{1}(A,\mun).$$
The universal exact sequence for $K$-theory with
coefficients yields a commutative diagram
$$\begin{array}{ccccccc}0\to & K_3(B^h)\otimes\Z/n&\to
&K_3(B^h;\Z/n)&\to &{}_nK_2(B^h)&\to 0\\
 &\downarrow & &\veq & &\downarrow & \\
0\to&K_3(A)\otimes\Z/n &\to&
K_3(A;\Z/n)&\by{}&{}_nK_2(A)&\to 0
\end{array}$$
Thus the right-most vertical arrow is surjective.
We then conclude from commutativity of the diagram:
$$\begin{array}{ccc}
\Het{1}(B^h,\mun) &\by{\simeq}&{}_nK_2(B^h)\\
 \veq  & &\quad\downarrow\mbox{\rm onto}\\
\Het{1}(A,\mun) &\by{\varphi}&{}_nK_2(A)\
\end{array}$$

\kern-14pt \E{proof}

\B{schol}\label{phi-bar} Let $A$ and $k$ be as in
Lemma~\ref{HoobH1}.  If $n$ is even, assume that $k$
contains a square root of $-1$.
If $\beta$ is a Bott element in $K_2(k;\Z/n)$ mapping to
$\zeta\in K_1(k)=k^*$, then
multiplication by $\beta$ lifts the map $\varphi$ to a map
$$\bar\varphi\colon\Het{1}(A,\mun)\to K_3(A;\Z/n).$$
This map is a split injection, because the \'etale Chern class
satisfies $c^{et}_2\bar\varphi=-1$.
\E{schol}
\B{proof} The assumption on $k$ implies that $\beta$ exists
and has order $n$, so $\bar\varphi(a)=\{ a,\beta\}$ is
well-defined and lifts $\varphi$.
The product formula (see \cite[Theorem 3.2(ii)]{W-chern}) states
that $c_2(\{a,\beta\})=-[a]\otimes\zeta$
in $A^*/A^{*n}\otimes\mu_n(k) \cong \Het{1}(\Spec(A),\mun)$
for every $a\in A^*$.
Thus up to sign $c_2$ is a left inverse of $\bar\varphi$.
\E{proof}


\B{thm}\label{nK2}
Let $X$ be a scheme of finite type over $\C$.
Then the \'etale Chern class defines
an isomorphism of Zariski sheaves:
$$c^{et}_2\colon\; {}_n\cK_2 \by{\simeq} \cH^1(\mun).$$
\E{thm}
\B{proof}  We saw in (2.4.2) that
$c^{et}_2\colon \cK_3(\Z/n) \to \cH^1(\mun)$ vanishes on $\cK_3$.
Hence $c^{et}_2$ is well-defined on ${}_n\cK_2$.
To verify that it is an isomorphism, we check the stalks at a
point $x\in X$.  If $A=\cO_{X,x}$ we see from \ref{phi-bar}
that the surjection $\varphi\colon \Het1(A,\mun)\to{}_nK_2(A)$
of  Lemma~\ref{HoobH1} satisfies $c^{et}_2\ \varphi=-1$.
Elementary algebra now implies that
$c^{et}_2$ is an isomorphism on ${}_nK_2(A)$
and hence on ${}_n\cK_2$.
 \E{proof}


\B{cor}\label{square} By (\ref{et-HD-sheaf}),
the following diagram commutes:
$$ 
\B{array}{ccc} {}_n\cK_2 &\longby\tau& \kern-4pt\cK_2 \\
c^{et}_2\downarrow\cong & &\downarrow{c_2}\\
\cH^1(\mun) &\longby{\delta}& \cH_\cD^2(2)
\E{array}$$ 
where $\tau$ is the obvious inclusion
and $\delta$ is defined in (\ref{triangle}) and (\ref{Rq}).
\E{cor}
%


\B{rmk} This gives the following explicit formula for $\delta$.
Given a unit $a\in A^*$, where $U=\Spec(A)$, write $[a]$
for the class of $c_1(a)$ in
$H_{an}^1(U_{\bul},\Z(1)_\cD)\cong
H_{an}^0(U_{\bul},\cO_{U_{\bul}}^*)$.
Then the product formula for $c_2$ shows that  $\delta$ sends
$a\otimes\zeta\in\cO_X^*(U)\otimes\mu_n\cong\cH^1(\mun)(U)$
to $[\zeta]\cup[a]\in H_{an}^2(U_{\bul},\Z(2)_\cD)=\cH_\cD^2(U)$.
\E{rmk}

\section{An exact sequence for $\cK_2$-cohomology}

We now give some exact sequences relating $H^1(X,\cK_2)$
and $H^2(X,\cK_2)$.  The first is a reinterpretation of
\cite[Theorem D]{PW3} in terms of hypercohomology.
Let $\cK_2^\bullet$ denote the complex $\cK_2\by{n}\cK_2$
concentrated in degrees 0 and 1.  The short exact sequence
$0\to\cK_2[-1]\to\cK_2^\bullet\to\cK_2\to0$
gives rise to a long exact sequence, reminiscient of
\cite[(4.4)]{S2}:
\B{equation}\B{array}{c}
0\to  H^0(X,\cK_2)/n \by\iota \H^1(X,\cK_2^\bullet)
\by{} H^1(X,\cK_2)\by{n}H^1(X,\cK_2)_{\mathstrut} \by\iota \\
\qquad\H^2(X,\cK_2^\bullet)\to H^2(X,\cK_2) \by{n}
H^2(X,\cK_2) \to\cdots^{\mathstrut}\E{array}
\E{equation}
 From this we extract short exact ``Kummer'' sequences, such as
\B{equation}\label{Kummer}
0\to H^1(X,\cK_2)/n \by\iota \H^2(X,\cK_2^\bullet)
\by{\pi} {}_nH^2(X,\cK_2) \to0.
\E{equation}
We also have the exact sequence of low degree terms in
the hypercohomology spectral sequence for $\cK_2^\bullet$,
the relevant part of which is:
\B{equation}\label{lowdeg}
H^0(X,\cK_2/n)\vlongby{d_2} H^2(X,{}_n\cK_2) \to
\H^2(X,\cK_2^\bullet)\by\eta H^1(X,\cK_2/n) \to H^3(X,{}_n\cK_2).
\E{equation}

\B{prop}\label{hyperK2}(\cite[Theorem D]{PW3})
Let $X$ be a quasi-projective over a field $k$ containing $\frac1n$.
Then:
\B{description} \item[{\it (a)}]  The $d_2$-differential
$H^0(X,\cK_2/n) \longby{d_2} H^2(X,{}_n\cK_2)$
in the hypercohomology spectral sequence (\ref{lowdeg})
is the composite
$$H^0(X,\cK_2/n) \by{\partial} H^1(X,n\cdot\cK_2)
\by{\partial} H^2(X,{}_n\cK_2)$$
of the boundary maps in the usual interlocking sequences
for $\cK_2$.
\item[{\it (b)}] If $X$ is a surface with isolated singularities,
the map $\pi$ in the Kummer sequence (\ref{Kummer}) for
$\H^2(X,\cK_2^\bullet)$ factors through the surjection $\eta$ in
the hypercohomology spectral sequence (\ref{lowdeg}).
\E{description}
\E{prop}
\B{proof} Part (a) is a special case of a more general result
which we have isolated in Lemma~\ref{hyper} below.
For part (b), it suffices to show that
the following diagram commutes.
$$\begin{array}{ccccccccc} 
H^0(X,\cK_2/n)\kern-5pt&
\by{\gamma}\kern-5pt&
H^2(X,{}_n\cK_2)\kern-8pt &\by{\beta}\kern-7pt&
H^1(X,\cK_2)/n\kern-7pt&\to\kern-7pt&H^1(X,\cK_2/n)\kern-5pt &
\to\kern-5pt& {}_nH^2(X,\cK_2) \kern-2pt\to0\\
\veq &
&\veq &&\iota\downarrow &&\veq &&\\
H^0(X,\cK_2/n)\kern-5pt&
\by{d_2}\kern-5pt&
H^2(X,{}_n\cK_2)\kern-8pt &\to\kern-7pt&
\H^2(X,\cK_2^\bullet)
\kern-7pt&\to\kern-7pt& H^1(X,\cK_2/n)
\kern-5pt&\to\kern-5pt&\ 0\hfill\\
\end{array} $$
The top row is the exact sequence of
\cite[Theorem D]{PW3}, the bottom row is the exact sequence of
low degree terms (\ref{lowdeg}) and the vertical arrow $\iota$
comes from the Kummer sequence (\ref{Kummer}).
Since $\cK_2\to\cK_2/n$ factors through $\cK^\bullet[1]$,
the right square commutes.  The left square commutes by part (a).
The map $\beta$ is constructed as follows.
Let $\cL^\bullet$ denote the subcomplex
$\cK_2\to n\cdot\cK_2$ of $\cK_2^\bullet$;
$\cL^\bullet$ is quasi-isomorphic to ${}_n\cK_2$.
The inclusion of
$n\cdot\cK_2[-1]$ into $\cL^\bullet$ induces a natural map
$H^1(X,n\cdot\cK_2)\to H^2(X,{}_n\cK_2)$, and
we know that this map is onto by \cite[4.8.1]{PW3}.
We showed in \cite[Proposition 4.9]{PW3} that
$H^1(X,n\cdot\cK_2)\to H^1(X,\cK_2)$ factors through
this surjection, and the induced map is $\beta$.
Thus $\iota\beta$ is induced from the composite map
$n\cdot\cK_2\to\cK_2\to\cK_2^\bullet[1]$
upon taking $H^1$.  But this composite map
factors through the subcomplex $\cL^\bullet[1]$
of $\cK^\bullet[1]$, so it follows that the left square commutes.
\E{proof}

Here is the general result about hypercohomology
spectral sequences used to prove part (a) above.
It works for any topos $X$.

\B{lemma}\label{hyper} For any sheaf $\cF$, let $\cC^\bullet$
denote the complex $\cF \by{n} \cF$ concentrated in degrees 0 and 1.
Then up to the sign $(-1)^{p}$, the $d_2$-differentials
$$H^p(X,\cF/n)=H^p(X,H^1\cC) \longrightarrow
  H^{p+2}(X,H^0\cC)=H^{p+2}(X,{}_n\cF)$$
in the hypercohomology spectral sequence of $\cC^\bullet$
are the composites
$$H^p(X,\cF/n)\by{\partial}H^{p+1}(X,n\cdot\cF)
	\by{\partial}H^{p+2}(X,{}_n\cF)$$
of the boundary maps $\partial$ associated
respectively to the exact sequences
$$0\to n\cdot\cF\to\cF\to\cF/n\to0,\qquad
	0\to{}_n\cF\to\cF\to n\cdot\cF\to0$$
\E{lemma}
\B{proof}
Given injective resolutions ${}_n\cF\to\cI^\bullet$,
$n\cdot\cF\to\cJ^\bullet$ and $\cF/n\cF\to\cK^\bullet$
we can form injective resolutions
$\cF\to\cE^{0\bullet}=\cI^\bullet\oplus\cJ^\bullet$ and
$\cF\to\cE^{1\bullet}=\cJ^\bullet\oplus\cK^\bullet$
using the Horseshoe Lemma.
These form the two columns of a Cartan-Eilenberg resolution
$\cE^{\bullet\bullet}$ of the complex $\cF\by{n}\cF$;
by the sign trick, the single horizontal differential
in this complex is $(-1)^{p}$ times the projection/inclusion
$I^p\oplus\cJ^p\to \cJ^p\to \cJ^p\oplus\cK^p$.

Given a class $[s]\in H^p(X,\cF/n)$, represent it by
$s\in H^0(X,K^p)$ with $\partial s=0$ in $H^0(X,K^{p+1})$.
Applying $\partial^v$ to $(0,s)\in H^0(X,J^p\oplus K^p)$
gives an element $(t,0)$ of $H^0(X,J^{p+1}\oplus K^{p+1})$.
Thus $\partial\colon H^p(X,\cF/n)\to H^{p+1}(X,n\cdot\cF)$
sends $[s]$ to $[t]$.  Applying $\partial^v$
to $(0,t)\in H^0(X,I^p\oplus J^p)$ gives $(u,0)$ for some
$u\in H^0(X,I^{p+1})$.  By construction, $u$ is a cycle in
$I^\bullet$ and
$\partial\colon H^{p+1}(X,n\cdot\cF)\to H^{p+2}(X,{}_n\cF)$
sends $[t]$ to $[u]$.

Now the hypercohomology spectral sequence arises from the row
filtration on the Cartan-Eilenberg resolution $\cE^{\bullet\bullet}$.
Since the pair $((0,(-1)^{p}t),(0,s))\in\mbox{Tot}^p(\cI)$
has $(((-1)^pu,0),(0,0))$ for its boundary,
the $d_2$-differential in the spectral sequence
takes $[s]$ to $(-1)^{p}[u]$.
\E{proof}


We are now going to connect Proposition~\ref{hyperK2}
with \'etale cohomology using $c^{et}_2$.  For this we need to
resort to some standard topological constructions.
Our main result will be Theorem~\ref{NH3} below.

Recall from (2.4.1) that there is a simplicial presheaf $K$
on $X_{zar}$ such that $\pi_qK(U)=K_{q}(U)$.
Let $\tilde K(U)$ be the universal covering space of the
basepoint component of $K(U)$; $\tilde K$ is a simplicial presheaf
by \cite[8.3 or 16.4]{May}.  Let $\tilde K^{(2)}(U)$ denote
the second layer of the Postnikov tower of $\tilde K(U)$,
defined in \cite[8.1]{May}; it is an Eilenberg-MacLane complex
of type $(K_2U,2)$ and $\tilde K^{(2)}$ is a simplicial presheaf.
Moreover by \cite[8.2]{May} there are Kan fibrations
$\tilde K^{(2)} \leftarrow \tilde K \to K$.

Now let $\tilde L(U)$ denote the homotopy fiber of the map
$\tilde K(U) \by{n} \tilde K(U)$,
and let $M(U)$ denote the homotopy fiber of the map
$\tilde K^{(2)}(U) \by{n} \tilde K^{(2)}(U)$.
Each $\tilde L(U)$ is a connected space with
$\pi_1\tilde L(U)=K_2(U)/n$
and $\pi_q\tilde L(U)=K_{q+1}(U;\Z/n)$ for $q\ge2$,
while $M(U)$ has only two nontrivial homotopy groups:
$\pi_1M(U)=K_2(U)/n$ and $\pi_2M(U)={}_nK_2(U)$.
In fact, it is not hard to see that $M(U)$ is homotopy
equivalent to the simplicial space obtained by applying the
Dold-Kan theorem to the presheaf of chain complexes
$K_2 \by{n} K_2$ concentrated in degrees 2 and 1.

We can perform the above constructions so that
there is a commutative diagram of simplicial presheaves
(in which the diagram (\ref{fib}) forms the right side):
\B{equation}\label{grid}\B{array}{ccccccc}
M &\leftarrow& \tilde L &\to& L &\vlongby{C^{ss}_2}& \Omega\cE \\
\downarrow && \downarrow && \downarrow &&\delta\downarrow\enskip\\
\tilde K^{(2)}&\leftarrow& \tilde K &\to& K &\longby{C_2^{ss}}&\cD\\
n\downarrow\enskip&&n\downarrow\enskip&&n\downarrow\enskip&&
					n\downarrow\enskip\\
\tilde K^{(2)}&\leftarrow& \tilde K &\to& K &\vlongby{C_2^{ss}}&\cD
\E{array}\E{equation}

\B{num}
Given any simplicial presheaf $F$ on $X$, the
{\it generalized sheaf cohomology groups} $\H^q(X,F)$
were defined for $q\le0$ by Brown and Gersten \cite[p.280]{BG}.
(The homotopy categories of simplicial presheaves and
simplicial sheaves are equivalent by \cite[2.8]{J}.
In particular, if $\tilde F$ is the simplicial sheaf
associated to $F$ then $\H^q(X,F)=\H^q(X,\tilde F)$.)

If $F$ is the simplicial sheaf associated by the Dold-Kan theorem
to a cochain complex $\cF$ (concentrated in negative degrees),
then $\H^q(X,F)\cong \H_{zar}^{q}(X,\cF)$ for $q\le0$
by \cite[p.281]{BG}.  Since the simplicial sheaf associated to
$\cK_2^\bullet[2]$ is the sheafification of $M$ we have
$\H^q(X,M)=\H_{zar}^{q+2}(X,\cK_2^\bullet)$.
Similarly, by (2.4.2) we have
$$\H^q(X,\Omega\cE)=\H^{q-1}(X,\cE)\cong \Het{3-q}(X,\mun).$$
In particular, diagram (\ref{grid}) induces maps
\B{equation}\label{branch}
\H^2(X,\cK_2^\bullet)=\H^0(X,M) \stackrel\lambda\leftarrow
\H^0(X,\tilde L)\longby{c_2^{et}} \H^0(X,\Omega\cE) = \Het3(X,\mun).
\E{equation}

If $F$ is a simplicial presheaf on $X$, we write $\tilde\pi_qF$
for the sheaf associated to the presheaf $U\mapsto \pi_qF(U)$.
For example, we have
$$\tilde\pi_qM = \cases{\cK_2/n& if $q=1$; \cr
		      {}_n\cK_2& if $q=2$; \cr
			0      & else.\cr}  \qquad
\tilde\pi_q\tilde L=\cases{\cK_2/n& if $q=1$; \cr
       	           \cK_{q+1}(\Z/n)& if $q\ge2$; \cr
			   0      & else.\cr}
$$

Now recall that $X$ is quasi-projective over $\C$.
By \cite[Theorem 3]{BG} there is a ``Brown-Gersten''
spectral sequence in the fourth quadrant:
$$E_2^{pq}=H_{zar}^p(X,\tilde\pi_{-q}F)
	\Rightarrow\H^{p+q}(X,F).$$
In general, this spectral sequence is ``fringed''
\cite[p.285]{BG}, but since all the $F$ we consider here
are infinite loop spaces this fringing does not affect
$\H^q(X,F)$ for $q\le0$.
\E{num}

\B{example}
Here is an example of the fringing phenomenon.
If $F$ is associated to a cochain complex $\cF$,
with $\cF^q=0$ for $q>0$, then
it is well known that the Brown-Gersten spectral sequence
for $F$ is the same as the hypercohomology spectral sequence
for $\cF$.  For example, the simplicial sheaf $\cE$ was defined
in (2.4.2) as being associated to $\tau^{\le0}\R\omega_*\Z/n[2i]$.
The hypercohomology spectral sequence of this complex
coincides with the Leray spectral sequence
for $\Het{2i+*}(X,\mu_n^{\otimes i})$ in the region $q\le0$.
Thus it is a fringed spectral sequence converging in the region
$p+q\le0$.  The line $p+q=+1$ converges to the kernel of
$\Het{2i+1}(X,\mun)\to H^0(X,\cH^{2i+1}(\mu_n^{\otimes i}))$.

On the other hand, the sheafification $\tilde M$ of $M$ is
associated to the complex of sheaves $\cK_2^\bullet[2]$.
Hence the Brown-Gersten spectral sequence for $M$
is the same as the hypercohomology spectral sequence for
$\cK_2^\bullet[2]$, and there is no fringing effect.
\E{example}

\medskip
Any morphism $E\to F$ of simplicial presheaves induces a morphism of
Brown-Gersten spectral sequences. 
Thus (\ref{branch}) gives us a commutative diagram:
\def\rod{\stackrel{c^{et}_2\quad\cong}{\kern-20pt
         \hbox to115pt{\rightarrowfill}\kern-35pt}}
\B{equation}\label{maze}\B{array}{ccccc}
H^0(X,\cK_2/n) && \rod        && H^0(X,\cH^2(\mun)) \\
   \veq       &&              &&     \veq     \\
H^0(X,\tilde\pi_1M) &\stackrel\cong\leftarrow&
	H^0(X,\tilde\pi_1\tilde L)
&\vlongby{c^{et}_2\enspace\cong}& H^0(X,\tilde\pi_2\cE) \\
d_2\downarrow\quad && d_2\downarrow\quad && d_2\downarrow\quad \\
H^2(X,\tilde\pi_2M) &\leftarrow& H^2(X,\tilde\pi_2\tilde L)
&\vlongby{c^{et}_2}& H^2(X,\tilde\pi_3\cE) \\
   \veq       &&              &&     \veq     \\
H^2(X,{}_n\cK_2)&& \rod       && H^2(X,\cH^1(\mun)) \\
\E{array}\E{equation}
(The bottom square of (\ref{maze}) commutes because,
as noted in (2.4.2), the Chern class map
$c_2^{et}\colon\cK_3(\Z/n)\to\cH^1(\mun)$
factors through ${}_n\cK_2$.)

The following description of the differential in the Leray
spectral sequence was suggested in \cite[(0.4)]{PW3}.

\B{prop}\label{diff}
If we identify $\cK_2/n$ with $\cH^2(\mun)$ by \ref{HoobH2}
and ${}_n\cK_2$ with $\cH^1(\mun)$ by \ref{nK2}, then the
differential $d_2\colon H^0(X,\cH^2(\mun)\to H^2(X,\cH^1(\mun)$
in the Leray spectral sequence for $\Het*(X,\mun)$ becomes
identified with the differential in \ref{hyperK2}(a), \ie
$$H^0(X,\cH^2(\mun))\cong H^0(X,\cK_2/n) \by{\partial}
H^1(X,n\cdot\cK_2)\by{\partial} H^2(X,{}_n\cK_2)
\cong H^2(X,\cH^1(\mun)).
$$
\E{prop}

\B{proof}
The left vertical map in (\ref{maze}) is the differential in the
hypercohomology spectral sequence for $\H^{*}(X,\cK_2^\bullet)$
by Example 7.4, and was described in Proposition~\ref{hyperK2}(a).
Again by Example~7.4, the right vertical map in (\ref{maze})
is the corresponding differential
in the Leray spectral sequence for $\Het{4+*}(X,\mun)$.
A diagram chase on (\ref{maze}),
starting at $H^0(X,\tilde\pi_1\tilde L)$, yields the result.
\E{proof}

\medskip

\B{defi} Following Suslin \cite{S2},
we define $N\Het3(X)$ to be the kernel
of the natural map
$\Het{3}(X,\mu_n^{\otimes 2})\to H^0(X,\cH^3(\mun))$.
Here $X$ can be any scheme in which $n$ is invertible.
Of course, when $X$ is a surface over an algebraically closed
field the sheaf $\cH^3(\mun)$ vanishes and we have
$N\Het3(X)=\Het3(X,\mun)$.
\E{defi}

\smallskip
The following result was proven by Suslin \cite[p. 19]{S2}
for smooth varieties.  It is a partial answer to
\cite[Question~2]{BV1} and was conjectured in \cite[(0.4)]{PW3}.

\B{thm} \label{NH3}
Let $X$ be a surface with isolated singularities over
a field $k$ containing an algebraically closed field and $\frac1n$.
Then
$$N\Het3(X)\cong\H^2(X,\cK_2\kern-3pt\by{n}\kern-3pt\cK_2).$$
In particular, by (\ref{Kummer})
there is a functorial short exact sequence:
$$0\to H^1(X,\cK_2)/n\to N\Het3(X)\to{}_nCH_0(X)\to 0.$$
\E{thm}

\B{proof}
Since $X$ is a surface, the Brown-Gersten spectral sequences
associated to the simplicial presheaves in (\ref{branch})
have only three nonzero columns. Using the computations
given in (7.3) for $\tilde\pi_qM$ and $\tilde\pi_q\tilde L$,
the resulting exact sequences form the rows of a commutative diagram.

\B{equation}\label{M-L-et}\B{array}{cccccccc}
\kern-8pt H^0(X,\cK_2/n)&\kern-3pt\by{d_2}&
\kern-3pt H^2(X,{}_n\cK_2) \kern-6pt&\to&
\kern-6pt \H^2(X,\cK_2^\bullet)
\kern-2pt&\to&\kern-6pt H^1(X,\cK_2/n) \kern-6pt&\to0\\
\kern-8pt\uparrow\cong&&\uparrow\mbox{\rm onto}&&
\uparrow\lambda&&\uparrow\cong&\\
\kern-8pt H^0(X,\cK_2/n)&\kern-3pt\by{d_2}&
\kern-3pt H^2(X,\cK_3(\Z/n)) \kern-6pt&\to&
\kern-6pt \H^0(X,\tilde L)
\kern-2pt&\to&\kern-6pt H^1(X,\cK_2/n) \kern-6pt&\to0\\
\kern-8pt\downarrow\cong&&
c_2^{et}\downarrow\mbox{\rm onto}\kern1.5em
&&\downarrow\kern1em&&\downarrow\cong&\\
\kern-8pt H^0(X,\cH^2(\mun))&\kern-3pt\by{d_2}&
\kern-3pt H^2(X,\cH^1(\mun)) \kern-6pt&\to&\kern-6pt N\Het3(X)
\kern-2pt&\to&\kern-6pt H^1(X,\cH^2(\mun)) \kern-6pt&\to0
\E{array}\E{equation}
The outside vertical maps are isomorphisms by \ref{HoobH2}.
The two vertical maps marked `onto' in (\ref{M-L-et}) are
actually split surjections with the same kernel,
and are identified by Lemma~\ref{HoobH1}
since $\varphi\colon\cH^1(\mun)\to{}_n\cK_2$ yields an
isomorphism on $H^2$.
Indeed, by \ref{phi-bar} we know that the map
$c_2^{et}\colon\cK_3(\Z/n)\to\cH^1(\mun)$ is a
surjection, split up to sign by
$\bar\varphi\colon \cH^1(\mun)\to\cK_3(\Z/n)$.
A diagram chase on (\ref{M-L-et}) shows that the two maps
$\H^0(X,\tilde L)\longby{\lambda}\H^2(X,\cK_2^\bullet)$
and $\H^0(X,\tilde L)\longby{c_2^{et}} N\Het3(X)$
are both onto with the same kernel.  Thus the quotients
$\H^2(X,\cK_2^\bullet)$ and $N\Het3(X)$ are isomorphic.
\E{proof}


%

\B{cor}\label{H3} If $k$ is algebraically closed
then the short exact sequence is:
$$0\to H^1(X,\cK_2)/n\to \Het{3}(X,\mu_n^{\otimes 2})
	\to{}_nCH_0(X)\to 0.$$
\E{cor}


\B{thm} \label{indivisible} Let $X$ be a normal projective
surface over an algebraically closed field $k$. Let $\ell$
a prime number, $\ell\neq$char$(k)$. Then
$$H^1(X,\cK_2)\otimes \Q_{\ell}/\Z_{\ell} =0\quad\mbox{and}
\quad \Het{3}(X,\Q_\ell/\Z_\ell) \cong CH_0(X)_{\ell-tors}$$
 \E{thm}
\B{proof} Choose a resolution of singularities
$\pi\colon X'\to X$.
Passing to the direct limit as $\nu\to\infty$, with
$n=\ell^\nu$, the short exact sequences of
Corollary~\ref{H3} become the rows of the commutative diagram
$$\begin{array}{ccccccc}
0\to & H^1(X,\cK_2)\otimes\Q_\ell/\Z_\ell
&\to&\Het{3}(X,\Q_\ell/\Z_\ell)&\to &CH_0(X)_{\ell-tors}&\to 0\\
 &\downarrow & &\downarrow & &\downarrow\cong & \\
0\to & H^1(X',\cK_2)\otimes\Q_\ell/\Z_\ell
&\to&\Het{3}(X',\Q_\ell/\Z_\ell)&\to &CH_0(X')_{\ell-tors}&\to 0.
\end{array}$$
The right-hand vertical map is an isomorphism by the
Collino-Levine Theorem \cite{C2} \cite{L-Alb}.
By \cite{CTR}, we have
 $H^1(X',\cK_2)\otimes\Q_{\ell}/\Z_{\ell} =0$.
Therefore it suffices to show that
$$\Het3(X,\Q_{\ell}/\Z_{\ell})\cong
\Het3(X',\Q_{\ell}/\Z_{\ell}).$$



There is a Mayer--Vietoris sequence for $\ell$-adic cohomology
similar to (\ref{exnorm}) for the square (\ref{birsquare}).
This yields an exact sequence
$$0\to T\to\Het{3}(X,\Z_{\ell})\to \Het{3}(X',\Z_{\ell})\to0$$
with $T = \Het{2}(E,\Z_{\ell})/\im(\Het{2}(X',\Z_{\ell}))$.
The proof of Proposition~\ref{pure} goes through in the
$\ell$-adic setting as well to show that $T$
 is a torsion group (\cf \cite[2.1]{C2}).
Since we also have
$\Het{4}(X,\Z_\ell)\cong \Het{4}(X',\Z)\cong\Z_\ell^c$,
the universal coefficient theorem yields the result:
$$\Het{3}(X,\Q_{\ell}/\Z_{\ell})\cong
\Het{3}(X,\Z_\ell)\otimes\Q_{\ell}/\Z_{\ell}\cong
\Het{3}(X',\Z_\ell)\otimes\Q_{\ell}/\Z_{\ell}\cong
\Het{3}(X',\Q_{\ell}/\Z_{\ell}).$$

\kern-24pt
\E{proof}


\section{Proof of the Main Theorem}

Let $X$ be a complex projective surface.
In Lemma~\ref{degree}
we constructed the Abel--Jacobi map
$\rho\colon A_0(X)\to J^2(X)$.
Our Main Theorem, stated in the Introduction,
states that $\rho$ induces an isomorphism
$A_0(X)_{tors}\cong J^2(X)_{tors}$.
We now proceed to prove the Main Theorem.

If $X$ is a normal surface then the result
$A_0(X)_{tors}\cong J^2(X)_{tors}$ is a paraphrase of
the theorem of Levine and Collino (see \cite{C2}, \cite{L-Alb})
that $A_0(X)_{tors}\cong J^2(\tilde X)_{tors}$ for any
resolution of singularities $\tilde X\to X$, because
$J^2(X)\cong J^2(\tilde X)$ by Corollary~\ref{abnorm}.

Granting the normal case, we shall establish the general case
of our Main Theorem by comparing a singular surface $X$ with its
normalization $\tilde X$.  For this, we need the following
crucial Lemma.
Let $\cH_{an}^2(\Z)$ denote the Zariski sheaf on $X$
associated to the presheaf $U\mapsto H_{an}^2(U,\Z)$

\B{lemma}\label{BVS}  Let $X$ be an irreducible proper
surface over $\C$.  Then the following composition is zero.
$$H_{\cD}^2(X,\Z(2))\by\varepsilon
H_{an}^2(X,\Z)\to H^0(X,\cH_{an}^2(\Z))$$
\E{lemma}
\B{proof} By Lemma~\ref{filt} the image of $\varepsilon$ is
the torsion subgroup of $H_{an}^2(X,\Z)$.
However, the sheaf $\cH_{an}^2(\Z)$ and hence its global sections
are torsion free  by \cite[Cor. 3]{BVS1}.
\E{proof}

\B{prop}\label{crux} If $X$ is an irreducible proper surface
over $\C$, the following natural map is zero.
$$H^0(X,\cK_2)\to H^0(X,\cK_2/n)$$
\E{prop}

\B{proof}  By Proposition~\ref{sheafHD} the natural map
$H_\cD^2(X,\Z(2))\to H^0_{zar}(X,\cH_\cD^2(2))$ is an isomorphism.
The Proposition follows from Lemma~\ref{BVS} and a chase on
the following diagram, the left part of which commutes
by (2.4.1).  
$$\B{array}{ccccccc}
K_2(X) &\by{c_2}& H_\cD^2(X,\Z(2))&\by\varepsilon&
H_{an}^2(X,\Z) &\to& H_{an}^2(X,\Z/n) \\
\downarrow&&\downarrow\cong &&\downarrow&&\downarrow\\
H^0(X,\cK_2) &\by{c_2}& H^0(X,\cH_\cD^2(2))&\to&
H^0(X,\cH_{an}^2(\Z)) &\to& H^0(X,\cH^2(Z/n)).
\E{array}$$
\vskip-17pt
\E{proof}

\B{rmk} When $X$ is a {\it smooth} proper variety over an
algebraically closed field of characteristic zero,
Proposition~\ref{crux} was proven by Colliot-Th\'el\`ene
and Raskind \cite{CTR},
and also by H. Esnault \cite{Esn} over $\C$.
\E{rmk}

\B{prop}\label{SK1}
Let $Z$ be a scheme which is proper over $\C$.
If $Z$ is either a curve or a normal surface then
\B{description}
\item[{\it i)}] $c_2\colon H^1(Z,\cK_2)_{tors}\cong H^2(Z,\Q/\Z)$
\item[{\it ii)}] $H^1(Z,\cK_2)\otimes \Q/\Z = 0$
\E{description}
\E{prop}
\goodbreak
\B{proof} The hypothesis on $Z$ allows us to use
\ref{HD-tors} for the isomorphism
$H^2(Z,\Q/\Z)\cong H_\cD^3(X,\Z(2))_{tors}$.
When $Z$ is a curve both assertions follow from
Theorem~\ref{sing} and this remark. When $Z$ is a normal surface,
part {\it ii)\,} was proven in Theorem~\ref{indivisible}.
In order to prove part {\it i)\,} for a normal surface $Z$,
we apply $H^1$ to Corollary~\ref{square}
and combine with the diagram of Corollary~\ref{morph}
to get a commutative diagram for each $n$:

\B{equation}\label{square2}\B{array}{ccc}
H^1(Z,{}_n\cK_2)&\by{\strut\tau_n} &H^1(Z,\cK_2)_{n-tors}\\
c_2^{et}\downarrow\cong & &\downarrow{c_2}\quad\\
H^1(Z,\cH^1(\mun))&\by\delta & H^1(Z,\cH_\cD^2(2))_{n-tors}\\
\downarrow && \downarrow\qquad \\
\Het2(Z,\mun) &\by\delta & H_\cD^3(Z,\Z(2))_{n-tors}.
 \E{array}\E{equation}
Taking the direct limit as $n\to\infty$ turns
$\mun$ into $\Q/\Z$.  Since
$H^2(Z,\Q/\Z)$ is the torsion subgroup of
$H_\cD^3(Z,\Z(2))$ by Proposition~\ref{HD-tors},
we have a commutative diagram

$$\B{array}{cccccccc}
&\kern-5pt H^1(Z,\cK_{2,tors})\kern-5pt
&\by{\tau}&\kern-3pt H^1(Z,\cK_2)_{tors}&
\kern-5pt\to & \kern-6pt\coker(\tau)&\to0&\\
&c_2^{et}\downarrow\;\cong&&\downarrow c_2&&\kern-5pt\downarrow&&\\
0\to\kern-3pt&\kern-5pt H^1(Z,\cH^1(\Q/\Z))\kern-5pt&\to&\kern-3pt
\Het2(Z,\Q/\Z)&\kern-5pt\to& \kern-6pt H^0(Z,\cH^2(\Q/\Z))
&\kern-7pt\longby{d_2}&\kern-8pt H^2(Z,\cH^1(\Q/\Z))\\
\E{array}$$
in which the bottom row is exact by Corollary~\ref{morph}.
Therefore in order to prove $(i)$
we are reduced to the claim that
$$\coker\tau
\cong \ker H^0(Z,\cH^2(\Q/\Z))\by{d_2} H^2(Z,\cH^1(\Q/\Z))$$

For each $n$, let $\gamma_n$ denote the composition
$H^0(Z,\cK_2/n)\kern-1.5pt\by\partial\kern-2pt H^1(Z,n\cdot\cK_2)
\kern-1.5pt\by\partial\kern-2pt H^2(Z,{}_n\cK_2)$
in the usual interlocking long exact sequences
\B{equation}\label{interlock}\B{array}{ccccccc}
H^1(Z,{}_n\cK_2)&\by{\mathstrut\tau_n}&H^1(Z,\cK_2)
&\to&H^1(Z,n\cdot\cK_2)&\by\partial&H^2(Z,{}_n\cK_2)\\
&&&&\veq&&\\
H^0(Z,\cK_2)&\by0& H^0(Z,\cK_2/n)
&\stackrel\partial\hookrightarrow&
H^1(Z,n\cdot\cK_2) &\to &H^1(Z,\cK_2).
\E{array}\E{equation}
The arrow marked `0' in this diagram is the zero map
by Proposition~\ref{crux}.
The other zig-zag composition in (\ref{interlock}),
from $H^1(Z,\cK_2)$ to 
$H^1(Z,\cK_2)$,
is multiplication by $n$. 
It follows from (\ref{interlock}) that
$$\ker(\gamma_n) \cong H^0(Z,\cK_2/n)\cap \im(H^1(Z,\cK_2))
\cong \frac{H^1(Z,\cK_2)_{n-tors}}{H^1(Z,{}_n\cK_2)}
= \coker(\tau_n).
$$
By Proposition \ref{hyperK2}(a), $\gamma_n$ is
the differential $d_2$ in the hypercohomology
spectral sequence for $\cK_2\by{n}\cK_2$.
By Proposition~\ref{diff}, we may also identify $\gamma_n$ with the
$d_2$-differential in the Leray spectral sequence for
$\Het{*}(Z,\mun)$.  Passing to the direct limit 
we obtain the claimed formula:
$\coker\tau = \lim\limits_{n\to\infty}\coker\tau_n\cong
\lim\limits_{n\to\infty} \ker\gamma_n=\ker(d_2)$.
\E{proof}


We are now ready to prove our Main Theorem for an arbitrary
projective surface $X$.  Letting $\tilde X$ be its normalization
and $Y$ a subscheme chosen as in Theorem \ref{M-V-K},
we have a Mayer-Vietoris Sequence in $K$-theory,
and also for Deligne cohomology by \ref{M-V-D}.
Taking the torsion subgroups of the diagram in
Corollary~\ref{SKtoD} yields the following
commutative diagram (in which we have abbreviated the
left-hand terms for legibility).

\B{equation}\label{diagram}\begin{array}{cccccccc}
\kern-5pt\biggl\{{SK_1(\tilde X)\oplus\atop SK_1(Y)}\biggr\}_{tors}
\kern-7pt&\to\kern-6pt & SK_1(\tilde Y)_{tors}
\kern-5pt&\to\kern-5pt&
SK_0(X)_{tors} \kern-5pt&\to\kern-7pt&
SK_0(\tilde X)_{tors} & \kern-5pt\to0\\
\kern-5pt\downarrow\mbox{}\kern-5pt
& &\downarrow\cong & &\downarrow\cong&&\downarrow\cong&\\
\kern-5pt\biggl\{{H^1(\tilde X)\oplus\atop H^1(Y,\cK_2)\ }
\biggr\}_{tors} \kern-7pt&\to\kern-6pt &
H^{1}(\tilde Y\kern-2pt,\cK_2)_{tors}
\kern-5pt&\to\kern-5pt&
H^{2}(X,\cK_2)_{tors} \kern-5pt&\to\kern-7pt&
H^{2}(\tilde X,\cK_2)_{tors} & \kern-5pt\to0\\
\kern-5ptc_2\downarrow\cong\ &&c_2\downarrow\cong
&&c_2\downarrow\ \ &&c_2\downarrow\cong \ &\\
\kern-5pt\biggl\{{H_\cD^3(\tilde X)\atop H_\cD^3(Y)}\biggr\}_{tors}
\kern-7pt&\to\kern-6pt & H_\cD^3(\tilde Y\kern-2pt,\Z(2))_{tors}
\kern-5pt&\to\kern-5pt&
H_\cD^4(X,\Z(2))_{tors}\kern-5pt&\to\kern-7pt&
H_\cD^4(\tilde X,\Z(2))_{tors}& \kern-5pt\to0
\end{array}
\E{equation}
Some discussion of diagram (\ref{diagram}) is in order.
The 3 isomorphisms between the terms in the top two rows
come from \ref{SKtoD}.
The 2 vertical maps in the lower left of (\ref{diagram})
are isomorphisms by Proposition \ref{SK1}.
The lower right vertical map
$H^{2}(\tilde X,\cK_2)_{tors}\cong J^2(\tilde X)_{tors}
\cong H_\cD^4(\tilde X,\Z(2))_{tors}$ is an isomorphism
because $\tilde X$ is normal.

The bottom row of (\ref{diagram}) is exact, because
by Proposition \ref{HD-tors} it is isomorphic to
$$H^2(\tilde X,\Q/\Z)\oplus H^2(Y,\Q/\Z)
\to H^2(\tilde Y,\Q/\Z) \to
H^3(X,\Q/\Z) \to H^3(\tilde X,\Q/\Z) \to0.
$$
The top two rows of (\ref{diagram}) are exact except
at $SK_1(\tilde Y)_{tors}$ and $H^1(\tilde Y,\cK_2)_{tors}$
by Proposition~\ref{SK1} and the
elementary lemma~\ref{tors-exact} below,
whose proof is left as an exercise.
The 5-lemma implies that we have an isomorphism
$$c_2\colon H^{2}(X,\cK_2)_{tors}\cong H_\cD^4(X,\Z(2))_{tors}$$
and this finishes the proof of our Main Theorem.
\hfil$\bullet$

\B{rmk}In order for the diagram chase of (\ref{diagram}) to work,
it suffices to know the crude surjectivity of the left vertical
map as $n\to\infty$:
$$H^1(\tilde X,\cK_2)_{tors}\oplus H^{1}(Y,\cK_2)_{tors}
\longby{c_2} H^2(\tilde X,\Q/\Z)\oplus H^2(Y,\Q/\Z).$$
\E{rmk}

\B{lemma}\label{tors-exact}
Let $A\to B\to C\to D$ be an exact sequence of abelian
groups.  If $A\otimes\Q/\Z=0$ then the following sequence is exact.
$$B_{tors} \to C_{tors} \to D_{tors}$$
\E{lemma}

Here is a motivic version of our Main Theorem.
For a 1-motive $M=(L,A,T,J,u)$ we let $M_{tors}$ denote the
extension of torsion subgroups.
$$0\to T_{tors}\to J_{tors}\to A_{tors}\to 0$$
Then our Main Theorem says that
$Alb(X)_{tors}$ can be described via {\it algebraic} zero-cycles,
\ie that $J^2(X)_{tors}$ is isomorphic to $A_0(X)_{tors}$
in a way compatible with normalization and desingularization.

\B{schol} Let $Alb(X)$ the Albanese 1--motive of
a projective surface. We then have the
following identification of $Alb(X)_{tors}$:
$$\begin{array}{ccccccc}
0\to&(\Q/\Z)^s &\to&A_0(X)_{tors} & \to &
A_0(\tilde X)_{tors}&\to 0\\
&\veq &&\downarrow\cong&&\downarrow\cong& \\
0\to&(\Q/\Z)^s &\to & J^2(X)_{tors}&
\to & J^2(\tilde X)_{tors} &\to 0. \end{array}$$
\E{schol}

\B{rmk} If $X$ is an {\it affine} surface over $\C$ then
$CH_0(X)=A_0(X)$ is uniquely divisible.  Indeed, the fact that
$A_0(X)_{tors} =0$ was proven in \cite[Theorem 2.6]{L2}.
And divisibility of $CH_0(X)=SK_0(X)$ is classical, probably
attributable to Murthy:
Every smooth point $x$ on $X$ is in the image of a map
$j\colon C\to X$ in which $C$ is a smooth affine curve.
The group $\Pic(C)$ is divisible, and the class of $x$
is in the image of the map $j_*\colon \Pic(C)\to SK_0(X)$.

Since $H^3(X,\C)=0$ as well, we also have $J^2(X)=0$.
Thus Roitman's Theorem holds by default in the affine case.
\E{rmk}


\vspace{2cm}

\B{flushright}{\sc DIMA - Universit\`a di Genova, Italy}\\[4pt]
{\sc Math. Dept. - Rutgers University, New Brunswick, USA}
\E{flushright}

\end{document}